\documentclass[reqno,11pt]{amsart}
\usepackage{amsmath}
\usepackage{amsxtra}
\usepackage{amscd}
\usepackage{amsthm}
\usepackage{amsfonts}
\usepackage{amssymb}
\usepackage{eucal}
\usepackage[matrix,arrow,curve]{xy}
\usepackage[dvipsnames]{xcolor}

\theoremstyle{plain}
\newtheorem{Thm}[subsection]{Theorem}
\newtheorem{Cor}[subsection]{Corollary}
\newtheorem{Lem}[subsection]{Lemma}
\newtheorem{Prop}[subsection]{Proposition}
\newtheorem{Conj}[subsection]{Conjecture}

\theoremstyle{definition}
\newtheorem{Def}[subsection]{Definition}

\theoremstyle{remark}

\newtheorem{Rem}[subsection]{Remark}

\newtheorem{conj}{Conjecture}

\numberwithin{equation}{section}
\renewcommand{\rm}{\normalshape}

\newif\ifShowLabels
\ShowLabelstrue
\newdimen\theight
\def\TeXref#1{%
    \leavevmode\vadjust{\setbox0=\hbox{{\tt
        \quad\quad  {\small \rm #1}}}%
    \theight=\ht0
    \advance\theight by \lineskip
    \kern -\theight \vbox to
    \theight{\rightline{\rlap{\box0}}%
    \vss}%
    }}%

\ShowLabelsfalse

\renewcommand{\sec}[2]{\section{#2}\label{S:#1}%
    \ifShowLabels \TeXref{{S:#1}} \fi}
\newcommand{\ssec}[2]{\subsection{#2}\label{SS:#1}%
    \ifShowLabels \TeXref{{SS:#1}} \fi}

\newcommand{\refs}[1]{Section ~\ref{S:#1}}
\newcommand{\refss}[1]{Section ~\ref{SS:#1}}

\newcommand{\reft}[1]{Theorem ~\ref{T:#1}}
\newcommand{\refl}[1]{Lemma ~\ref{L:#1}}
\newcommand{\refp}[1]{Proposition ~\ref{P:#1}}

\newcommand{\refd}[1]{Definition ~\ref{D:#1}}

\newcommand{\refe}[1]{\eqref{E:#1}}

\newenvironment{thm}[1]%
    { \begin{Thm} \label{T:#1}  \ifShowLabels \TeXref{T:#1} \fi }%
    { \end{Thm} }

\renewcommand{\th}[1]{\begin{thm}{#1} \sl }
\renewcommand{\eth}{\end{thm} }

\newenvironment{lemma}[1]%
    { \begin{Lem} \label{L:#1}  \ifShowLabels \TeXref{L:#1} \fi }%
    { \end{Lem} }
\newcommand{\lem}[1]{\begin{lemma}{#1} \sl}
\newcommand{\elem}{\end{lemma}}

\newenvironment{propos}[1]%
    { \begin{Prop} \label{P:#1}  \ifShowLabels \TeXref{P:#1} \fi }%
    { \end{Prop} }
\newcommand{\prop}[1]{\begin{propos}{#1}\sl }
\newcommand{\eprop}{\end{propos}}

\newenvironment{corol}[1]%
    { \begin{Cor} \label{C:#1}  \ifShowLabels \TeXref{C:#1} \fi }%
    { \end{Cor} }
\newcommand{\cor}[1]{\begin{corol}{#1} \sl }
\newcommand{\ecor}{\end{corol}}

\newenvironment{defeni}[1]%
    { \begin{Def} \label{D:#1}  \ifShowLabels \TeXref{D:#1} \fi }%
    { \end{Def} }
\newcommand{\defe}[1]{\begin{defeni}{#1} \sl }
\newcommand{\edefe}{\end{defeni}}

\newenvironment{remark}[1]%
    { \begin{Rem} \label{R:#1}  \ifShowLabels \TeXref{R:#1} \fi }%
    { \end{Rem} }
\newcommand{\rem}[1]{\begin{remark}{#1}}
\newcommand{\erem}{\end{remark}}

\newenvironment{conjec}[1]%
    { \begin{Conj} \label{Co:#1}  \ifShowLabels \TeXref{Co:#1} \fi }%
    { \end{Conj} }
\renewcommand{\conj}[1]{\begin{conjec}{#1} \sl }
\newcommand{\econj}{\end{conjec}}

\newcommand{\eq}[1]%
    { \ifShowLabels \TeXref{E:#1} \fi
       \begin{equation} \label{E:#1} }
\newcommand{\eeq}{ \end{equation} }

\newcommand{\prf}{ \begin{proof} }
\newcommand{\epr}{ \end{proof} }

\newcommand{\nc}{\newcommand}

\newcommand{\iso}{\stackrel{\sim}{\longrightarrow}}
\nc{\HC}{{\mathcal{HC}}}
\nc{\on}{\operatorname}
\nc{\BA}{{\mathbb{A}}}
\nc{\BC}{{\mathbb{C}}}
\nc{\BG}{{\mathbb{G}}}
\nc{\BM}{{\mathbb{M}}}
\nc{\BN}{{\mathbb{N}}}
\nc{\BQ}{{\mathbb{Q}}}
\nc{\BP}{{\mathbb{P}}}
\nc{\BR}{{\mathbb{R}}}
\nc{\BZ}{{\mathbb{Z}}}
\nc{\BS}{{\mathbb{S}}}

\nc{\CA}{{\mathcal{A}}}
\nc{\CB}{{\mathcal{B}}}
\nc{\CalC}{{\mathcal C}}
\nc{\CalD}{{\mathcal D}}
\nc{\CE}{{\mathcal{E}}}
\nc{\CF}{{\mathcal{F}}}
\nc{\CG}{{\mathcal{G}}}
\nc{\CH}{{\mathcal{H}}}
\nc{\CK}{{\mathcal{K}}}
\nc{\CL}{{\mathcal{L}}}
\nc{\CM}{{\mathcal{M}}}
\nc{\CMM}{{\mathcal{M}^{\operatorname{gen}}_\hbar(-\rho)}}
\nc{\CN}{{\mathcal{N}}}
\nc{\CO}{{\mathcal{O}}}
\nc{\CP}{{\mathcal{P}}}
\nc{\CQ}{{\mathcal{Q}}}
\nc{\CR}{{\mathcal{R}}}
\nc{\CS}{{\mathcal{S}}}
\nc{\CT}{{\mathcal{T}}}
\nc{\CU}{{\mathcal{U}}}
\nc{\CV}{{\mathcal{V}}}
\nc{\CW}{{\mathcal{W}}}
\nc{\CX}{{\mathcal{X}}}
\nc{\CY}{{\mathcal{Y}}}
\nc{\CZ}{{\mathcal{Z}}}

\nc{\gen}{{\operatorname{gen}}}
\nc{\cM}{{\check{\mathcal M}}{}}
\nc{\csM}{{\check{\mathcal A}}{}}
\nc{\obM}{{\overset{\circ}{\mathbf M}}{}}
\nc{\oCA}{{\overset{\circ}{\mathcal A}}{}}
\nc{\obA}{{\overset{\circ}{\mathbf A}}{}}
\nc{\ooM}{{\overset{\circ}{M}}{}}
\nc{\osM}{{\overset{\circ}{\mathsf M}}{}}
\nc{\vM}{{\overset{\bullet}{\mathcal M}}{}}
\nc{\nM}{{\underset{\bullet}{\mathcal M}}{}}
\nc{\obD}{{\overset{\circ}{\mathbf D}}{}}
\nc{\cp}{{\overset{\circ}{\mathbf p}}{}}
\nc{\ofZ}{{\overset{\circ}{\mathfrak Z}}{}}
\nc{\oCZ}{{\overset{\circ}{\mathcal Z}}{}}

\nc{\fa}{{\mathfrak{a}}}
\nc{\fb}{{\mathfrak{b}}}
\nc{\ff}{{\mathfrak{f}}}
\nc{\fg}{{\mathfrak{g}}}
\nc{\fgl}{{\mathfrak{gl}}}
\nc{\fh}{{\mathfrak{h}}}
\nc{\fj}{{\mathfrak{j}}}
\nc{\fm}{{\mathfrak{m}}}
\nc{\fn}{{\mathfrak{n}}}
\nc{\fu}{{\mathfrak{u}}}
\nc{\fp}{{\mathfrak{p}}}
\nc{\frr}{{\mathfrak{r}}}
\nc{\fs}{{\mathfrak{s}}}
\nc{\ft}{{\mathfrak{t}}}
\nc{\fy}{{\mathfrak{y}}}
\nc{\ofT}{{\overline{\mathfrak T}}}
\nc{\ofS}{{\overline{\mathfrak S}}}
\nc{\fsl}{{\mathfrak{sl}}}
\nc{\hsl}{{\widehat{\mathfrak{sl}}}}
\nc{\hgl}{{\widehat{\mathfrak{gl}}}}
\nc{\hg}{{\widehat{\mathfrak{g}}}}
\nc{\chg}{{\widehat{\mathfrak{g}}}{}^\vee}
\nc{\hn}{{\widehat{\mathfrak{n}}}}
\nc{\chn}{{\widehat{\mathfrak{n}}}{}^\vee}

\nc{\fA}{{\mathfrak{A}}}
\nc{\fB}{{\mathfrak{B}}}
\nc{\fC}{{\mathfrak{C}}}
\nc{\fD}{{\mathfrak{D}}}
\nc{\fE}{{\mathfrak{E}}}
\nc{\fF}{{\mathfrak{F}}}
\nc{\fG}{{\mathfrak{G}}}
\nc{\fI}{{\mathfrak{I}}}
\nc{\fJ}{{\mathfrak{J}}}
\nc{\fK}{{\mathfrak{K}}}
\nc{\fL}{{\mathfrak{L}}}
\nc{\fM}{{\mathfrak{M}}}
\nc{\fN}{{\mathfrak{N}}}
\nc{\frP}{{\mathfrak{P}}}
\nc{\fQ}{{\mathfrak Q}}
\nc{\fR}{{\mathfrak R}}
\nc{\fS}{{\mathfrak S}}
\nc{\fT}{{\mathfrak{T}}}
\nc{\fU}{{\mathfrak{U}}}
\nc{\fV}{{\mathfrak{V}}}
\nc{\fW}{{\mathfrak{W}}}
\nc{\fZ}{{\mathfrak{Z}}}

\nc{\bb}{{\mathbf{b}}}
\nc{\bc}{{\mathbf{c}}}
\nc{\be}{{\mathbf{e}}}
\nc{\bj}{{\mathbf{j}}}
\nc{\bn}{{\mathbf{n}}}
\nc{\bp}{{\mathbf{p}}}
\nc{\bq}{{\mathbf{q}}}
\nc{\bs}{{\mathbf{s}}}
\nc{\bt}{{\mathbf{t}}}
\nc{\bv}{{\mathbf{v}}}
\nc{\bx}{{\mathbf{x}}}
\nc{\by}{{\mathbf{y}}}
\nc{\bw}{{\mathbf{w}}}
\nc{\bA}{{\mathbf{A}}}
\nc{\bB}{{\mathbf{B}}}
\nc{\bC}{{\mathbf{C}}}
\nc{\bK}{{\mathbf{K}}}
\nc{\bD}{{\mathbf{D}}}
\nc{\bH}{{\mathbf{H}}}
\nc{\bI}{{\mathbf{I}}}
\nc{\bM}{{\mathbf{M}}}
\nc{\bN}{{\mathbf{N}}}
\nc{\bO}{{\mathbf{O}}}
\nc{\bQ}{{\mathbf Q}}
\nc{\bS}{{\mathbf{S}}}
\nc{\bT}{{\mathbf{T}}}
\nc{\bV}{{\mathbf{V}}}
\nc{\bW}{{\mathbf{W}}}
\nc{\bX}{{\mathbf{X}}}
\nc{\bP}{{\mathbf{P}}}
\nc{\bZ}{{\mathbf{Z}}}

\nc{\sa}{{\mathsf{a}}}
\nc{\sfb}{{\mathsf{b}}}
\nc{\sA}{{\mathsf{A}}}
\nc{\sB}{{\mathsf{B}}}
\nc{\sC}{{\mathsf{C}}}
\nc{\sD}{{\mathsf{D}}}
\nc{\sF}{{\mathsf{F}}}
\nc{\sK}{{\mathsf{K}}}
\nc{\sM}{{\mathsf{M}}}
\nc{\sO}{{\mathsf{O}}}
\nc{\sQ}{{\mathsf{Q}}}
\nc{\sP}{{\mathsf{P}}}
\nc{\sR}{{\mathsf{R}}}
\nc{\sT}{{\mathsf{T}}}
\nc{\sV}{{\mathsf{V}}}
\nc{\sW}{{\mathsf{W}}}
\nc{\sX}{{\mathsf{X}}}
\nc{\sZ}{{\mathsf{Z}}}
\nc{\sfp}{{\mathsf{p}}}
\nc{\sr}{{\mathsf{r}}}
\nc{\st}{{\mathsf{t}}}
\nc{\sv}{{\mathsf{v}}}
\nc{\sy}{{\mathsf{y}}}
\nc{\sfc}{{\mathsf{c}}}
\nc{\sd}{{\mathsf{d}}}
\nc{\sg}{{\mathsf{g}}}
\nc{\sfl}{{\mathsf{l}}}

\nc{\BK}{{\bar{K}}}

\nc{\tA}{{\widetilde{\mathbf{A}}}}
\nc{\tB}{{\widetilde{\mathcal{B}}}}
\nc{\tg}{{\widetilde{\mathfrak{g}}}}
\nc{\tG}{{\widetilde{G}}}
\nc{\TM}{{\widetilde{\mathbb{M}}}{}}
\nc{\tO}{{\widetilde{\mathsf{O}}}{}}
\nc{\tU}{{\widetilde{\mathfrak{U}}}{}}
\nc{\TZ}{{\tilde{Z}}}
\nc{\tZ}{\widetilde{Z}{}}
\nc{\tx}{{\tilde{x}}}
\nc{\tbv}{{\tilde{\bv}}}
\nc{\tfP}{{\widetilde{\mathfrak{P}}}{}}
\nc{\tz}{{\tilde{\zeta}}}
\nc{\tmu}{{\tilde{\mu}}}

\nc{\td}{\ddot{\underline{d}}{}}
\nc{\tzeta}{\widetilde{\zeta}{}}
\nc{\hd}{{\widehat{\underline{d}}}}
\nc{\hG}{{\widehat{G}}}
\nc{\hBP}{\widehat{\mathbb P}{}}
\nc{\hQ}{{\widehat{Q}}}
\nc{\hsM}{\widehat{\mathsf M}{}}
\nc{\hfM}{\widehat{\mathfrak M}{}}
\nc{\hCP}{\widehat{\mathcal P}{}}
\nc{\hCR}{\widehat{\mathcal R}{}}
\nc{\hCS}{{\widehat{\mathcal S}}}
\nc{\hfZ}{\widehat{\mathfrak Z}{}}

\nc{\urho}{\underline{\rho}}
\nc{\uB}{\underline{B}}
\nc{\uC}{{\underline{\mathbb{C}}}}
\nc{\ui}{\underline{i}}
\nc{\ofP}{{\overline{\mathfrak{P}}}}
\nc{\oZ}{\overset{\circ}{Z}{}}
\nc{\oA}{\overset{\circ}{\mathbb A}{}}

\nc{\hrho}{{\hat{\rho}}}

\nc{\unl}{\underline}
\nc{\ol}{\overline}
\nc{\one}{{\mathbf{1}}}
\nc{\two}{{\mathbf{t}}}

\nc{\Tot}{{\mathop{\operatorname{\rm Tot}}}}
\nc{\Hilb}{{\mathop{\operatorname{\rm Hilb}}}}
\nc{\End}{{\operatorname{End}}}
\nc{\Ext}{{\operatorname{Ext}}}
\nc{\Hom}{{\operatorname{Hom}}}
\nc{\Sym}{{\operatorname{Sym}}}
\nc{\CHom}{{\mathop{\operatorname{{\mathcal{H}}om}}}}
\nc{\defi}{{\mathop{\operatorname{\rm def}}}}
\nc{\length}{{\mathop{\operatorname{\rm length}}}}

\nc{\Cliff}{{\mathsf{Cliff}}}
\nc{\Fib}{{\mathsf{Fib}}}
\nc{\Coh}{{\mathsf{Coh}}}
\nc{\FCoh}{{\mathsf{FCoh}}}

\nc{\reg}{{\text{\rm reg}}}

\nc{\cplus}{{\mathbf{C}_+}}
\nc{\cminus}{{\mathbf{C}_-}}
\nc{\cthree}{{\mathbf{C}_*}}
\nc{\Qbar}{{\bar{Q}}}
\nc{\Fl}{{{\mathcal F}\ell}}
\nc{\bh}{{\bar{h}}}
\nc{\bOmega}{{\overline{\Omega}}}
\nc\tGr{\widetilde{\Gr}}
\nc{\ul}{\underline}
\nc{\seq}[1]{\stackrel{#1}{\sim}}
\nc\ogu{\overline{G/U}}
\nc\chlam{\check{\lam}}

\nc\St{\operatorname{St}}
\nc{\bLambda}{{\boldsymbol{\Lambda}}}
\nc\uS{\underline{S}}
\nc\QM{\mathcal{QM}}
\nc{\chmu}{\check{\mu}}
\nc{\CHH}{{\CH\!\!\CH}}
\nc{\Lam}{\Lambda}
\nc{\calB}{\mathcal B}
\nc{\CC}{\mathbb C}
\nc{\ZZ}{\mathbb Z}
\nc{\alp}{\alpha}
\nc{\PP}{\mathbb P}
\nc\x{\times}
\nc\lam{\lambda}
\nc\gam{\gamma}
\nc\calF{\mathcal F}
\nc\disj{\operatorname{disj}}
\newcommand\calO{\mathcal O}
\newcommand\calL{\mathcal L}
\newcommand\la{\langle}
\newcommand\ra{\rangle}

\begin{document}
Dedicated to the memory of Andrei Zelevinsky

\bigskip

\author{Alexander Braverman, Galyna Dobrovolska, and Michael Finkelberg}
\title{Gaiotto-Witten superpotential and Whittaker D-modules on monopoles}

\begin{abstract}
  Let $G$ be an almost simple simply connected group over $\CC$. For a positive element $\alp$ of the coroot lattice of
  $G$ let $\oZ^{\alp}$ denote the space of maps from $\PP^1$ to the flag variety $\calB$ of $G$ sending $\infty\in \PP^1$
  to a fixed point in $\calB$ of degree $\alp$. This space is known to be isomorphic to the space of framed $G$-monopoles on
$\BR^3$ with maximal symmetry breaking at infinity of charge $\alp$.

  In~\cite{fkmm} a system of (\'etale, rational) coordinates on $\oZ^{\alp}$ is introduced. In this note we compute various known structures
  on $\oZ^{\alp}$ in terms of the above coordinates. As a byproduct we give a natural interpretation of the Gaiotto-Witten
  superpotential studied in \cite{gw} and relate it to the theory of Whittaker D-modules discussed in \cite{g}.
\end{abstract}

\maketitle 

\sec{Intro}{Introduction}

\ssec{}{Zastava spaces}
Let $G$ be an almost simple simply connected algebraic group over $\CC$. We denote by $\calB$ the flag variety of
$G$. Let us also fix a pair of opposite Borel subgroups $B$, $B_-$ whose intersection is a maximal torus $T$ (thus we have
$\calB=G/B=G/B_-$).

Let $\Lam$ denote the cocharacter lattice of $T$; since $G$ is assumed to be simply connected, this is also the coroot lattice of $G$.
We denote by $\Lam_+\subset \Lam$ the sub-semigroup spanned by positive coroots. We say that $\alp\geq \beta$ (for $\alp,\beta\in \Lam$)
if $\alp-\beta\in\Lam_+$.

It is well-known that $H_2(\calB,\ZZ)=\Lam$ and that an element $\alp\in H_2(\calB,\ZZ)$ is representable by an algebraic curve
if and only if $\alp\in \Lam_+$.
Let $\oZ^{\alp}$ denote the space of maps $\PP^1\to \calB$ of degree $\alp$ sending $\infty\in \PP^1$ to $B\in \calB$.
It is known~\cite{fkmm} that this is a smooth symplectic affine algebraic variety, which can be identified with the space of
framed $G$-monopoles on $\BR^3$ with maximal symmetry breaking at infinity of charge $\alp$~\cite{j},~\cite{j'}.

The scheme $\oZ^{\alp}$ is endowed with a number of remarkable structures
(listed below). On the other hand in~\cite{fkmm} the authors introduce a
system of (birational, \'etale) coordinates on $\oZ^{\alp}$. The purpose of
the present note is to compute how these structures look
like in the above coordinates. In particular, it turns out that the
Gaiotto-Witten superpotential~\cite{gw} admits a natural interpretation in
terms of Whittaker $D$-modules of~\cite{g}.

\ssec{}{Quasi-maps}
The scheme $\oZ^{\alp}$ has a natural partial compactification $Z^{\alp}$.
It can be realized as the
space of based {\em quasi-maps} of degree $\alp$; set-theoretically it can be described in the following way:
$$
Z^{\alp}=\bigsqcup\limits_{0\leq\beta\leq \alp} \oZ^{\beta}\x \BA^{\alpha-\beta},
$$
where for $\gam\in \Lam_+$ we denote by $\BA^{\gam}$ the space of all colored divisors $\sum\gam_i x_i$ with
$x_i\in \BA^1$, $\gam_i\in \Lam_+$ such that $\sum \gam_i=\gam$.
\ssec{symmetric}{A ``symmetric" definition of the Zastava space}
Fix
$\lam,\mu\in \Lam$.
Let us denote by $\oZ^{\lam,\mu}$ the scheme classifying the following data:

1) A $G$-bundle $\calF$ on $\PP^1$ with a trivialization at $\infty\in \PP^1$.

2) A $B$-structure $\calF_B$ on $\calF$ such that the induced $T$-bundle $\calF_{T,+}$ is of degree $\lam$. We require that $\CF_B$ is equal to
$B$ at $\infty$.

3) A $B_-$-structure $\calF_{B_-}$ on $\calF$ such that the induced $T$-bundle $\calF_{T,-}$ is of degree $\mu$. We require that $\CF_{B_-}$ is equal to
$B_-$ at $\infty$.

\noindent
It is easy to see that this is indeed a scheme. Moreover, we claim that $\oZ^{\lam,\mu}$ is naturally isomorphic to $\oZ^{\lam-\mu}$~(\refs{gaits}).

\ssec{}{Structures on the Zastava space}
It is easy to see that the
space $\oZ^{\alp}$ is endowed with the following structures (the precise
constructions are given in the main body of the paper):
\begin{enumerate}
  \item
  The scheme $\oZ^{\alp}$ possesses a natural symplectic structure~\cite{fkmm}.
  \item
  There is a natural morphism $\pi_{\alp}:\oZ^{\alp}\to \BA^{\alp}$. Moreover, given $\beta,\gam\in \Lam_+$ let
  $(\BA^{\beta}\x\BA^{\gam})_{\disj}$ denote the space of pairs of colored divisors of degrees $\beta$ and
  $\gam$ which are mutually disjoint. If $\alp=\beta+\gam$ then we have a natural \'etale map $(\BA^{\beta}\x\BA^{\gam})_{\disj}\to \BA^{\alp}$.
  The {\em factorization} is a
  canonical isomorphism
  $$
  \ff_{\beta,\gamma}:\
  (\BA^\beta\times\BA^\gamma)_{\on{disj}}\times_{\BA^\alpha}Z^\alpha\iso
  (\BA^\beta\times\BA^\gamma)_{\on{disj}}\times_{\BA^\beta\times\BA^\gamma}
  (Z^\beta\times Z^\gamma).
  $$
  We shall refer to the latter as the {\em factorization property} of Zastava.
  \item
  The Cartan involution on $G$ (which interchanges $B$ and $B_-$ and induces the map $t\mapsto t^{-1}$ on $T$) induces an involution $\iota$ on $\oZ^{\alp}$ (this is clear from the point of view of the definition of $\oZ^{\alp}$ given in \refss{symmetric}).
  \item
  Let $\partial Z^{\alp}=Z^{\alp}\backslash \oZ^{\alp}$.
Then $\partial Z^{\alp}$ is a Cartier divisor and moreover it is the divisor
  of zeros of some function $F_{\alp}$ on $Z^{\alp}$ which is invertible on $\oZ^{\alp}$ (this function is unique up to a multiplicative scalar).
  \item
  Fix $\lam,\mu\in \Lam$ such that $\lam-\mu=\alp$. Then for every simple root $\check\alp_i$ of $G$  we have canonical maps
  $\fE_{\lam,+,i}^{\alp}:\oZ^{\alp}\to H^1(\PP^1,\calO(\la-\check\alp_i, \lam\ra)$, $\fE_{\mu,-,i}^{\alp}:\oZ^{\alp}\to H^1(\PP^1,\calO(\la-\check\alp_i, \mu\ra)$.
  The precise definition is given in~\refs{gaits}, so let us just explain the definition for
  $G=SL(2)$ here. In this case $Z^{\lam,\mu}\simeq Z^{\alp}$ just classifies rank 2 vector bundles $\calF$ on $\PP^1$ with trivialized
  determinant together with two short exact sequences $0\to \calL_+\to \calF\to \calL_+^{-1}\to 0$ and $0\to \calL_-\to \calF\to \calL_-^{-1}\to 0$
  with $\deg \calL_+=-\lam, \deg\calL_-=-\mu$, where we identify the lattice $\Lam$ with $\ZZ$ in a natural way.
  In addition $\calF$ is endowed with a trivialization at $\infty$, which is compatible with $\calL_+$ and $\calL_-$; in particular
  $\calL_+$ and $\calL_-$ also get a trivialization at $\infty$ which allows us to identify them canonically with $\calO(-\lam)$ and
  $\calO(-\mu)$ (here we use a notation $\CO(n),\ n\in\BZ$, for a line bundle
on $\BP^1$ trivialized at $\infty\in\BP^1$).
Hence the above short exact sequences define elements in
$H^1(\PP^1,\calO(-2\lam))$ and $H^1(\PP^1,\calO(-2\mu))$.

 Let $\chi^{\lam}_{i,+}:\oZ^{\lambda,\mu}\x H^0(\PP^1,\calO(\la\lam,\check\alp_i\ra-2))\to \CC$ be the composition of
  $\fE_{\lambda,+,i}^\alpha$ and the natural pairing
$H^0(\PP^1,\calO(\la\lam,\check\alp_i\ra-2))\x H^1(\PP^1,\calO(-\la\lam,\check\alp_i\ra))\to
  \CC$. Note that an element of $H^0(\PP^1,\calO(\la\lam,\check\alp_i\ra-2)$ can be regarded as a polynomial
  $K_i$ of one variable $z$ of degree $\leq\la\lam,\check\alp_i\ra-2$.
Similarly, we let $\chi^{\mu}_{i,-}:\oZ^{\lambda,\mu}\x H^0(\PP^1,\calO(\la\mu,\check\alp_i\ra-2))\to \CC$ be the corresponding
  function (obtained by replacing $\fE_{\lambda,+,i}^\alpha$ with
$\fE_{\mu,-,i}^\alpha$).
We set $\fE_{\lam,+}^{\alp}$ to be the direct sum of all the $\fE_{\lam,+,i}^{\alp}$ and similarly for $\fE_{\mu,-}^{\alp}$
  (sometimes we shall drop the indices $\lam,\mu$ and $\alp$ when it does not lead to a confusion). Obviously the maps $\fE_+$ and
  $\fE_-$ are interchanged by the involution $\iota$.
\end{enumerate}

\ssec{}{Coordinates on Zastava}
A system of \'etale birational coordinates on $\oZ^{\alp}$ is introduced
in~\refss{igrek}. Let us recall the definition
for $G=SL(2)$. In this case $\oZ^{\alp}$ consists of all maps $\PP^1\to \PP^1$ of degree $\alp$ which send
$\infty$ to $0$. We can represent such a map by a rational function $\frac{R}{Q}$ where $Q$ is a monic polynomial of degree
$\alp$ and $R$ is a polynomial of degree $<\alp$. Let $w_1,\ldots,w_{\alp}$ be the zeros of $Q$. Set $y_r=R(w_r)$.
Then the functions $(y_1,\ldots,y_\alp,w_1,\ldots, w_{\alp})$ form a system of \'etale birational coordinates on $\oZ^{\alp}$.

For general $G$ the definition of the above coordinates is quite similar. In this case given a point in $\oZ^{\alp}$
 we can define polynomials $R_i,Q_i$ where $i$ runs through the set of vertices of the Dynkin diagram of $G$ and

(1) $Q_i$ is a monic polynomial of degree $\la\alp,\check\omega_i\ra$

(2) $R_i$ is a polynomial of degree $<\la\alp,\check\omega_i\ra$.

Hence, we can define (\'etale, birational) coordinates
$(y_{i,r},w_{i,r})$ where $i$ is as above and
$r=1,\ldots,\la\alp,\check\omega_i\ra$. It will be convenient for us to use slightly modified coordinates
$\by_{i,r}:=y_{i,r}\prod_{j\ne i}Q_j^{\langle\alpha_j,\check\alpha_i\rangle/2}(w_{i,r})$.
Then the main result of this note is the following
\th{main}
\begin{enumerate}
  \item
  The Poisson brackets of the modified coordinates (with respect to the symplectic structure defined in~\cite{fkmm}) are as follows:
  $$
  \{w_{i,r},w_{j,s}\}=0,\
  \{w_{i,r},\by_{j,s}\}=\check{d}_i\delta_{ij}\delta_{rs}\by_{j,s},\
  \{\by_{i,r},\by_{j,s}\}=0.
  $$
  \item
  (Recall that the boundary equation $F_\alpha$ is defined up to a
multiplicative constant.) We have $F_\alpha=\prod_{i,r}\by_{i,r}^{d_i}=
  \prod_{i,r}y_{i,r}^{d_i}\prod_{j\ne i}Q_j^{\alpha_j\cdot\alpha_i/2}(w_{i,r})$.
  \item
  Let us introduce yet another modified system of rational \'etale coordinates on $\oZ^\alpha$: we define
  \eq{etan}
  \fy_{i,r}:=\frac{y_{i,r}}
  {Q'_i(w_{i,r})}.
  \end{equation}
where $Q'_i$ stands for the derivative of the polynomial $Q_i(z)$.
  Then we have
  \eq{factor}
  \ff_{\beta,\gamma}(w_{i,r},\fy_{i,r})_{i\in I}^{1\leq r\leq a_i}=\\
  \left((w_{i,r},\fy_{i,r})_{i\in I}^{1\leq r\leq b_i},
  (w_{i,r},\fy_{i,r})_{i\in I}^{b_i+1\leq r\leq a_i}\right).
  \end{equation}
  \item
  The involution $\iota$ sends
  $(w_{i,r},\by_{i,r})$ to $(w_{i,r},\by_{i,r}^{-1})$.
  \item
  We have
  \eq{ext+}
  \begin{aligned}
  \chi^{\lam}_{i,+}(\unl{w},\unl{y},\unl{z})=\sum_{r=1}^{a_i}y_{i,r}^{-1}
  \frac{\prod_{j\ne i}Q_j^{-\langle\alpha_j,\check\alpha_i\rangle}(w_{i,r})}
  {Q'_i(w_{i,r})}K_i(w_{i,r})=\\
  \sum_{r=1}^{a_i}\by_{i,r}^{-1}
  \frac{\prod_{j\ne i}Q_j^{-\langle\alpha_j,\check\alpha_i\rangle/2}(w_{i,r})}
  {Q'_i(w_{i,r})}K_i(w_{i,r}).
  \end{aligned}
  \end{equation}
  Similarly,
  \eq{ext-}
  \chi^{\mu}_{i,-}(\unl{w},\unl{y},\unl{z})=
  \sum_{r=1}^{a_i}\by_{i,r}
  \frac{\prod_{j\ne i}Q_j^{-\langle\alpha_j,\check\alpha_i\rangle/2}(w_{i,r})}
  {Q'_i(w_{i,r})}K_i(w_{i,r}).
  \end{equation}
  \end{enumerate}
\eth

\rem{star}
The set of irreducible components $\on{Irr}^\alpha$ of the central factorization
fiber $\pi_\alpha^{-1}(\alpha\cdot0)\subset Z^\alpha$ is in a natural bijection
with the weight $\alpha$ component of the Kashiwara crystal
$\sB_{\check\fg}(\infty)$,~\cite[Section~14]{bfg}. The involution induced by
$\iota$ on $\bigsqcup_\alpha\on{Irr}^\alpha$ is nothing but the involution
$*:\ \sB_{\check\fg}(\infty)\to\sB_{\check\fg}(\infty)$ of~\cite[8.3]{ka}.
\erem

\ssec{relgwg}{Relation with the works of Gaiotto-Witten and Gaitsgory}
We keep the notation from \reft{main}. Let us observe that a monic polynomial $K(z)$ of degree $d$
is the same as a point in $\BA^{(d)}$. Thus if all $K_i$ are monic, together they form an point in
$\BA^{\lam-2\rho}$. Thus we may regard $\chi^{\lam}_{\pm}:=\sum_{i\in I}\chi^\lambda_{i,\pm}$ as functions on $\oZ^{\alp}\x \BA^{\lam-2\rho}$.

Let $\Lam=(\lam_1,\ldots,\lam_n)$ be an unordered collection of dominant
coweights whose sum is equal to $\lam-2\rho$.
Then $\Lam$ defines a locally closed subvariety $\oA^{\Lambda}$ in
$\BA^{\lam-2\rho}$ (namely, the moduli space of configurations of {\em distinct}
colored points $\unl{z}=(z_1,\ldots,z_n)$ so that the color of $z_i$ is
$\lambda_i$) and we denote by
$\chi^{\Lam}_{\pm}$ the restriction of $\chi^{\lam}_{\pm}$ to $\oZ^{\alp}\x \oA^{\Lam}$.
We now define the (multivalued) superpotentials
$\CW^{\Lam,\alp}_{\pm}:\fh^{\vee}\times \oZ^{\alp}\times \oA^{\Lam}\to \BA^1$ by setting
\eq{superpm}
\CW^{\Lam,\alp}_{\pm}=\sum_{1\leq n\leq N}\langle\lambda_n,h^*\rangle z_n-
\sum_{(i,r)}\langle\alpha_i,h^*\rangle w_{i,r}\pm\log F_{\alp}+\chi^{\Lam}_{\pm}+\sum_{1\leq m<n\leq N}\lambda_m\cdot\lambda_n\log(z_m-z_n).
\end{equation}
Note that all the summands except the 3rd and the 4th are pulled back from $\oA^{\alp}\times \BA^{\Lam}$.
Also, it is clear from the above definition that the exponential
of $\CW^{\Lam}_+$ is well defined as
a regular function on $\fh^\vee\times \oZ^{\alp}\times \oA^{\Lam}$.
In addition the involution $\iota$ transforms $\CW^{\Lam,\alp}_+$ to $\CW^{\Lam, \alp}_-$.

Let us now assume that $G=SL(2)$. Then it follows from
\reft{main} that the function $\CW^{\Lam,\alp}_-$ is exactly the Gaiotto-Witten superpotential studied in~\cite{gw}.
We shall from now on use this name for $\CW^{\Lam,\alp}_-$ for any $G$. We see that the exponential of
the Gaiotto-Witten super-potential is well-defined on $\fh^\vee\times \oZ^{\alp}\times \oA^{\Lam}$ (this is not immediately clear from the
coordinate description).

On the other hand, let $\varkappa\in \CC$ be an irrational number.
Then the work of Gaitsgory \cite{g} easily implies the following result:
\th{gaits-to-witten}
Let $M^{\varkappa,\alp,\Lam}_-$ denote  the $D$-module on $\oZ^{\alp}\times \oA^{\Lam}$
generated by the function $\exp(\varkappa\CW^{\Lam,\alp}_-)$.
Let $\pi^{\alp,\Lam}:\ \fh^\vee\times\oZ^{\alp}\times \oA^{\Lam}\to
\fh^\vee\times\BA^{\alp}\times \oA^{\Lam}$ be the corresponding morphism.
Then we have $\pi^{\alp,\Lam}_!(M^{\varkappa,\alp,\Lam}_-)=
\pi^{\alp,\Lam}_*(M^{\varkappa,\alp,\Lam}_-)$ and it is isomorphic to the
the minimal extension of the $D$-module on the open stratum generated by
the function
$$\prod_{1\leq n\leq N}\exp(\langle\lambda_n,\varkappa h^*\rangle z_n)\times
\prod_{(i,r)}\exp(-\langle\alpha_i,\varkappa h^*\rangle w_{i,r})\times$$
$$\times\prod_{(i,r)\ne(j,s)}(w_{i,r}-w_{j,s})^{\varkappa\alpha_i\cdot\alpha_j/2}\times
\prod_{(i,r),1\leq n\leq N}(z_n-w_{i,r})^{-\varkappa\alpha_i\cdot\lambda_n}\times
\prod_{1\leq m<n\leq N}(z_m-z_n)^{\varkappa\lambda_m\cdot\lambda_n}.
$$
\eth
\ssec{}{Remark} The above Theorem is essentially due to Gaiotto and Witten when restricted to the open stratum (in this case
it is not difficult to deduce it from the coordinate description of the
superpotential).
Interpreting the superpotential in terms of \refe{superpm} allows one to extend this statement to all of $\BA^{\alp}\times \oA^{\Lam}$
using the work of Gaitsgory. It would be interesting to find an interpretation of this refined statement in terms of
the Landau-Ginzburg model studied by Gaiotto and Witten.

\ssec{ack}{Acknowledgments} We are grateful to R.~Bezrukavnikov,
B.~Feigin, D.~Gaiotto,
D.~Gaitsgory, M.~Gekhtman, S.~Oblezin, L.~Rybnikov, V.~Schechtman and
A.~Uteshev for the useful discussions. A.B. was partially supported by the NSF and
by the Simons Foundation. The financial support from the Government of the Russian Federation within the framework of the implementation of the 5-100 Programme Roadmap of the National Research University  Higher School of Economics, AG Laboratory  is acknowledged by M.F.

\sec{Ex1}{Recollections about zastava}

\ssec{nota}{Notations}
$G$ is an almost simple simplyconnected complex algebraic Lie group.
We fix its Cartan and Borel subalgebras $T\subset B\subset G$ with the Lie
algebras $\fh\subset\fb\subset\fg$. The set of simple roots is denoted $I$;
the simple roots (resp. coroots) are denoted $\check\alpha_i$
(resp. $\alpha_i$), $i\in I$.
We fix a Weyl group invariant symmetric bilinear form $?\cdot?$ on
the Cartan Lie algebra $\fh$ such that the square length of a {\em short}
coroot is $\alpha_i\cdot\alpha_i=2$.
This bilinear form gives rise to an isomorphism $\fh^\vee\iso\fh$ so that
the root lattice $X$ generated by $\{\check\alpha_i\}_{i\in I}$ embeds into
$\fh$. We then have $\check\alpha_i\cdot\check\alpha_i\in\{2,1,\frac{2}{3}\}$,
and $\alpha_i\cdot\alpha_i\in\{2,4,6\}$.
We set $d_i=\frac{\alpha_i\cdot\alpha_i}{2}$. Let $d$ be the ratio of the square
lengths of the long and short coroots, so that $d\in\{1,2,3\}$. We set
$\check{d}_i=d/d_i$. Then $\langle\alpha_i,\check\alpha_j\rangle=
\frac{\alpha_i\cdot\alpha_j}{d_j}=d_i\check\alpha_i\cdot\check\alpha_j=
d\frac{\check\alpha_i\cdot\check\alpha_j}{\check{d}_i}$.

For $\alpha=\sum_{i\in I}a_i\alpha_i,\ a_i\in\BN$, we consider the corresponding
zastava space $Z^\alpha$ (see e.g.~\cite{bf}) with an open smooth subvariety
$Z^\alpha\supset\oZ^\alpha$: the moduli space of degree $\alpha$ based maps from
$C=\BP^1$ to the flag variety $\CB=G/B$ (also known as the moduli space of
framed $G$-monopoles on $\BR^3$ of topological charge $\alpha$ with the maximal
symmetry breaking at infinity). The complementary boundary divisor
is denoted $\partial Z^\alpha:=Z^\alpha\setminus\oZ^\alpha$.

\ssec{igrek}{Coordinates on zastava}
Let $z$ be a coordinate on $C=\BP^1$. We think of the zastava space $Z^\alpha$
in its Pl\"ucker embedding as of collections of degree
$\langle\alpha,\check\lambda\rangle$ $V_{\check\lambda}$-valued
polynomials (here $\check\lambda$ is a dominant weight, and $V_{\check\lambda}$
is the corresponding irreducible representation) such that the highest weight
component is of the form $z^{\langle\alpha,\check\lambda\rangle}+\ldots$ (the smaller
powers of $z$), and all the other weight components are of degree strictly
smaller than $\langle\alpha,\check\lambda\rangle$.
In particular, if $\check\lambda=\check\omega_i$, a fundamental weight,
then the highest weight component is denoted $Q_i$ (a monic polynomial of
degree $a_i=\langle\alpha,\check\omega_i\rangle$), and the prehighest weight
($=\check\omega_i-\check\alpha_i$) component is denoted $R_i$ (a polynomial
of degree $<a_i$). The polynomial $Q_i$ is determined uniquely by the
(unordered) set of its roots $w_{i,r},\ 1\leq r\leq a_i$. The ramified cover
$\varpi:\ \widehat{Z}{}^\alpha\to Z^\alpha$ is formed by all the orderings of the
roots of all the polynomials $Q_i,\ i\in I$. We have regular functions
$y_{i,r}:=R_i(w_{i,r})$ on $\widehat{Z}{}^\alpha$. According
to~\cite[Remark~2]{fkmm}, on the open subset where all the roots $w_{i,r},\
i\in I$ are distinct (and $\varpi$ is unramified), $\{w_{i,r},y_{i,r}\}$ form
a coordinate system (an open embedding into $\BA^{\langle\alpha,2\check\rho\rangle}$).

\ssec{modi}{A symplectic form and modified coordinates}
The main result of~\cite{fkmm} is a construction of a symplectic form on
$\overset{\circ}{Z}{}^\alpha$ which extends as a Poisson structure to $Z^\alpha$.
According to~\cite[Proposition~2]{fkmm}, the Poisson brackets of the coordinates
of~\refss{igrek} are as follows: $\{w_{i,r},w_{j,s}\}=0,\
\{w_{i,r},y_{j,s}\}=\check{d}_i\delta_{ij}\delta_{rs}y_{j,s},\ \{y_{i,r},y_{j,s}\}=
d\check\alpha_i\cdot\check\alpha_j\frac{y_{i,r}y_{j,s}}{w_{i,r}-w_{j,s}}$ for
$i\ne j$, and finally $\{y_{i,r},y_{i,s}\}=0$.

Following the private communications of S.~Oblezin and L.~Rybnikov, we
consider the modified rational \'etale coordinates
$\by_{i,r}:=y_{i,r}\prod_{j\ne i}Q_j^{\langle\alpha_j,\check\alpha_i\rangle/2}(w_{i,r})$
(they are regular only on the open subset where all the roots $w_{i,r},\
i\in I$ are distinct).

\lem{rybobl}
The Poisson brackets of the modified coordinates are as follows:
$\{w_{i,r},w_{j,s}\}=0,\
\{w_{i,r},\by_{j,s}\}=\check{d}_i\delta_{ij}\delta_{rs}\by_{j,s},\
\{\by_{i,r},\by_{j,s}\}=0$.
\elem

\prf
Straightforward.
\epr
Note that this is exactly the statement of \reft{main}(1).
\defe{loga}
We define the logarithmic coordinates $\sy_{i,r}:=\log\by_{i,r}$ on an
appropriate $\BZ^{|\alpha|}$-cover of the open subset of $\widehat{Z}{}^\alpha$
where all the roots $w_{i,r},\ i\in I$ are distinct, and $y_{i,r}\ne0$.
\edefe

\ssec{version}{A version of zastava}
Given $\lambda,\mu\in X_*(T)$, we consider the moduli stack $\oZ^{\lambda,\mu}$
of the following data: (a) a $G$-bundle $\CF_G$ on $C$ trivialized at
$\infty\in C$; (b) a reduction of $\CF_G$ to a $B$-bundle (a $B$-structure
on $\CF_G$) such that the induced $T$-bundle has degree $\lambda$,
and the fiber of the $B$-structure at $\infty\in C$ is $B\subset G$;
(c) a reduction of $\CF_G$ to a $B_-$-bundle (a $B_-$-structure
on $\CF_G$) such that the induced $T$-bundle has degree $\mu$,
and the fiber of the $B_-$-structure at $\infty\in C$ is $B_-\subset G$.

According to~\cite[Section~2]{bfgm}, $\oZ^{\lambda,\mu}$ is representable by
a scheme. More precisely, $\alpha:=\lambda-\mu$ is automatically a nonnegative
combination of positive coroots, and $\oZ^{\lambda,\mu}$ is isomorphic to the
zastava scheme $\oZ^\alpha$.

The Cartan involution of $G$ interchanging $B$ and $B_-$ and acting on
$T$ as $t\mapsto t^{-1}$ induces an isomorphism
$\iota:\ \oZ^{\lambda,\mu}\iso\oZ^{-\mu,-\lambda}$. The composition
$\oZ^\alpha\simeq\oZ^{\lambda,\mu}\stackrel{\iota}{\longrightarrow}\oZ^{-\mu,-\lambda}
\simeq\oZ^\alpha$ is a well defined involution $\iota:\ \oZ^\alpha\iso\oZ^\alpha$
(independent of the choice of a presentation $\alpha=\lambda-\mu$: the
independence is clear from the description of the identification $\oZ^{\lambda,\mu}\iso\oZ^\alpha$ of~\cite[Section~2]{bfgm}).

\sec{factoriza}{Factorization (Proof of \reft{main}(3))}

\ssec{facto}{Factorization in coordinates}
Recall the fundamental {\em factorization} property of zastava spaces.
For $\alpha=\beta+\gamma$ we have a natural morphism
$a:\ \BA^\beta\times\BA^\gamma\to\BA^\alpha$. An open subset
$(\BA^\beta\times\BA^\gamma)_{\on{disj}}\subset(\BA^\beta\times\BA^\gamma)$
is formed by the pairs $(D_\beta,D_\gamma)$
of disjoint divisors $D_\beta,D_\gamma\in\BA^1$. The {\em factorization} is a
canonical isomorphism $$\ff_{\beta,\gamma}:\
(\BA^\beta\times\BA^\gamma)_{\on{disj}}\times_{\BA^\alpha}Z^\alpha\iso
(\BA^\beta\times\BA^\gamma)_{\on{disj}}\times_{\BA^\beta\times\BA^\gamma}
(Z^\beta\times Z^\gamma).$$
We introduce yet another modified system of rational \'etale coordinates
on $Z^\alpha$: we define
\eq{eta}
\fy_{i,r}:=\frac{y_{i,r}}
{Q'_i(w_{i,r})}.
\end{equation}

Let $\beta=\sum_{i\in I}b_i\alpha_i,\ \gamma=\sum_{i\in I}c_i\alpha_i$, so that
$a_i=b_i+c_i$.

\prop{fact coor}
$\ff_{\beta,\gamma}(w_{i,r},\fy_{i,r})_{i\in I}^{1\leq r\leq a_i}=
\left((w_{i,r},\fy_{i,r})_{i\in I}^{1\leq r\leq b_i},
(w_{i,r},\fy_{i,r})_{i\in I}^{b_i+1\leq r\leq a_i}\right)$.
\eprop

\prf
We recall the construction of the factorization isomorphism.
Let $U$ stand for the unipotent radical of the Borel $B$, and let $U_-$ be
the unipotent radical of the opposite Borel (with the same Cartan torus $T$)
$B_-$. Let $\ol{G/U}$ stand for the affinization of the base affine space.
The quotient stack $U_-\backslash\ol{G/U}/T$ has an open dense point; and
the complement is a Cartier (Schubert) divisor $D$. Now $Z^\alpha$ is the moduli
space of degree $\alpha$ maps $C\to U_-\backslash\ol{G/U}/T$ (i.e. such that
the induced $T$-bundle on $C$ has degree $\alpha$) such that $\infty\in C$
goes to the complement of the Schubert divisor, see e.g.~\cite{bfgm}.

For $\phi\in Z^\alpha$, the pullback of the Schubert divisor $\phi^*D$ is
nothing but $\pi_\alpha(\phi)\in(C\setminus\infty)^\alpha$. Given
$\phi_1\in Z^\beta,\ \phi_2\in Z^\gamma$ with disjoint $\pi_\beta(\phi_1),
\pi_\gamma(\phi_2)$, we construct the corresponding $\phi\in Z^\alpha$ as
follows. Note that the disjointness condition guarantees that
$\CU_1:=C\setminus\phi_1^*D$ and $\CU_2:=C\setminus\phi_2^*D$ cover $C$, and
$\phi_1|_{\CU_1\cap\CU_2}=\phi_2|_{\CU_1\cap\CU_2}$ (the constant map to the point).
So we define $\phi$ by gluing $\phi_1$ and $\phi_2$ over $\CU_1\cap\CU_2$.

Now let us replace $G,U,U_-,T$ by $SL_2^i,U^i,U_-^i,T^i$ corresponding to the
$i$-th root. Then $\ol{SL_2^i/U^i}$ is isomorphic to a 2-dimensional vector
space $V_i$; the right action of $T^i$ is isomorphic to the scalar action of
$\BC^*$; the left action of $U_-^i$ is isomorphic to the one coming from the
natural left action of $SL_2^i$. We have the canonical homomorphisms
$\chi_i:\ U_-\twoheadrightarrow U_-^i$, and $\check\alpha_i:\ T\to T^i$.
We also have a natural projection
$\on{pr}_i:\ \overline{G/U}\twoheadrightarrow\ol{SL_2^i/U^i}$. In effect,
$\overline{G/U}$ in Pl\"ucker realization consists of collections of vectors
in the irreducible $G$-modules. In particular, each collection contains a
vector $v_{\check\omega_i}\in V_{\check\omega_i}$. So we set
$\on{pr}_i(v_{\check\lambda})_{\check\lambda\in X^*(T)^+}:=\on{pr}_i(v_{\check\omega_i})\in
V_{\check\omega_i}^{\on{Rad}P_i}=V_i$. It is straightforward to check that
$\on{pr}_i$ is $\chi_i:\ U_-\twoheadrightarrow U_-^i$-equivariant, and
$\check\alpha_i:\ T\to T^i$-equivariant. In other words, we have a morphism
of stacks
$\on{pr}_i:\ U_-\backslash\ol{G/U}/T\to U_-^i\backslash\ol{SL_2^i/U^i}/T^i$,
and the inverse image of the Schubert divisor
$D_i\subset U_-^i\backslash\ol{SL_2^i/U^i}/T^i$ lies inside the Schubert divisor
$D\subset U_-\backslash\ol{G/U}/T$ (in fact, this inverse image coincides with
the corresponding irreducible component of $D$). Hence we obtain the same
named projection $\on{pr}_i:\ Z^\alpha_\fg\to Z_{\mathfrak{sl}_2^i}^{a_i}$, and the
following diagram commutes:
\eq{Cd}
\begin{CD}
Z^\alpha_\fg @>\on{pr}_i>> Z_{\mathfrak{sl}_2^i}^{a_i}\\
@V\pi_\alpha VV      @V\pi_{a_i} VV\\
\BA^\alpha @>\on{pr}_i>> \BA^{(a_i)}
\end{CD}
\end{equation}
Moreover, the following diagram commutes as well:
\eq{cD}
\begin{CD}
(\BA^\beta\times\BA^\gamma)_{\on{disj}}\times_{\BA^\alpha}Z^\alpha_\fg @>\ff_{\beta,\gamma}>>
(\BA^\beta\times\BA^\gamma)_{\on{disj}}\times_{\BA^\beta\times\BA^\gamma}
(Z^\beta_\fg\times Z^\gamma_\fg)\\
@V\on{pr}_i VV      @V\on{pr}_i VV\\
(\BA^{(b_i)}\times\BA^{(c_i)})_{\on{disj}}\times_{\BA^{(a_i)}}Z^{a_i}_{\mathfrak{sl}_2^i}
@>\ff_{b_i,c_i}>>
(\BA^{(b_i)}\times\BA^{(c_i)})_{\on{disj}}\times_{\BA^{(b_i)}\times\BA^{(c_i)}}
(Z^{b_i}_{\mathfrak{sl}_2^i}\times Z^{c_i}_{\mathfrak{sl}_2^i})
\end{CD}
\end{equation}
Hence the proposition is reduced to the case of $\fg={\mathfrak{sl}_2^i}$
that will be dealt with in the next section.

\ssec{sltwo}{Factorization for $SL_2$}
In this section $G=SL_2^i$, and to unburden the notations we will write
$G,U,U_-,T$ for $SL_2^i,U^i,U_-^i,T^i$. We will use another point of view on
the factorization. Namely, we will think of
$Z^a\ni\phi:\ C\to U_-\backslash\ol{G/U}/T$ as of a $G$-bundle $\CF$ on
$C$ with a generalized $B$-structure, and a $U_-$-structure transversal to
the $B$-structure at $\infty\in C$. These generically transversal structures
define a generic trivialization of $\CF$, i.e. a point of the Beilinson-Drinfeld
Grassmannian $\on{Gr}_{BD}$. Moreover, since any $U_-$-bundle over $C$ is
trivial, $\CF$ is trivial too, and its trivialization at $\infty\in C$ extends
to a canonical global trivialization. Thus the above trivialization (coming
from two transversal structures) may be viewed as a rational function
$C\to G$; more precisely, as a rational function $f:\ C\to U_-$ (because of
the reduction to $U_-$) sending $\infty\in C$ to the neutral element of $U_-$.
Now recall that $G=SL_2^i$, and $U_-={\mathbb G}_a=\BA^1$. Then in the elementary
terms $f$ is nothing but $\frac{R_i}{Q_i}$.

Back to factorization, it arises from the factorization of the
Beilinson-Drinfeld Grassmannian. Given $G$-bundles $\CF_1,\CF_2$ with
trivializations $\sigma_1,\sigma_2$ defined on the open subsets
$\CU_1,\CU_2\subset C$ such that $\CU_1\cup\CU_2=C$ we construct a new bundle
$\CF$ with trivialization $\sigma$ on $\CU=\CU_1\cap\CU_2$ by gluing
$\CF_1|_{\CU_2}$ and $\CF_2|_{\CU_1}$ over $\CU$ where they are both trivialized.

Given $Z^b\ni\phi_1$ (resp. $Z^c\ni\phi_2$) corresponding to
$(\CF_1,\CU_1,\sigma_1)$ (resp. $(\CF_2,\CU_2,\sigma_2)$) and
$f_1=\frac{R_1}{Q_1}$ (resp. $f_2=\frac{R_2}{Q_2}$) we want to compute the
result of gluing $Z^a\ni\phi$ corresponding to $(\CF,\CU,\sigma)$ and
$f=\frac{R}{Q}$. Note that by the construction, the principal part of
$f$ at $C\setminus\CU_1$ (resp. $C\setminus\CU_2$) coincides with the
principal part of $f_1$ at $C\setminus\CU_1$ (resp. with that of $f_2$ at
$C\setminus\CU_2$). On the other hand, the rational function $f$ of degree $a$
vanishing at $\infty\in C$ is uniquely determined by its principal parts
at $(C\setminus\CU_1)\cup(C\setminus\CU_2)$. We conclude $f=f_1+f_2$.
This is equivalent to the desired formula of~\refe{eta} and~\refp{fact coor}
(since the principal part of $f$ at $w_{i,r}$, i.e. the residue of $fdz$ at
$w_{i,r}$, is given by the formula~\refe{eta}).

This completes the proof of the proposition.
\epr

\ssec{fact rev}{Another factorization}
Recall from~\refss{sltwo} that the factorization isomorphism
$$\ff_{\beta,\gamma}:\
(\BA^\beta\times\BA^\gamma)_{\on{disj}}\times_{\BA^\alpha}\oZ^\alpha\iso
(\BA^\beta\times\BA^\gamma)_{\on{disj}}\times_{\BA^\beta\times\BA^\gamma}
(\oZ^\beta\times\oZ^\gamma)$$ (\refss{facto}) is induced by the embedding
$\oZ^\alpha\hookrightarrow\on{Gr}_{BD}(U_-)\hookrightarrow\on{Gr}_{BD}(G)$.
Given $\unl{x}=\sum_m\alpha_m\cdot x_m\in\BA^\alpha$ the fiber
$\pi_\alpha^{-1}(\unl{x})$ goes under this embedding to
$\prod_m(\fT_0\cap\fS_{\alpha_m})\subset\on{Gr}_{BD}(G)$.
Here $\fT_0\subset\on{Gr}_{G,x_m}$ (resp. $\fS_{\alpha_m}\subset\on{Gr}_{G,x_m}$)
is the semiinfinite orbit $U_-(\CK_{x_m})\cdot0$
(resp. $U(\CK_{x_m})\cdot\alpha_m$), and $\CK_{x_m}\supset\CO_{x_m}$ is the local
field (resp. ring) around the point $x_m\in C$, and $\alpha_m\in\on{Gr}_G$
is a $T$-fixed point. Note that $\fT_0\subset\on{Gr}_G$ is canonically
isomorphic to $\on{Gr}_{U_-}\subset\on{Gr}_G$.

We also have a natural embedding $\oZ^{0,-\alpha}\hookrightarrow\on{Gr}_{BD}(G)$
sending the fiber over $\unl{x}$ to
$\prod_m(\fT_{-\alpha_m}\cap\fS_0)\subset\on{Gr}_{BD}(G)$. Note that
$\fS_0\subset\on{Gr}_G$ is canonically isomorphic to $\on{Gr}_U\subset\on{Gr}_G$.
Under the identification $\oZ^\alpha\simeq\oZ^{0,-\alpha}$ the factorization
of $\on{Gr}_U$ induces the factorization
$$\ff^+_{\beta,\gamma}:\
(\BA^\beta\times\BA^\gamma)_{\on{disj}}\times_{\BA^\alpha}\oZ^\alpha\iso
(\BA^\beta\times\BA^\gamma)_{\on{disj}}\times_{\BA^\beta\times\BA^\gamma}
(\oZ^\beta\times\oZ^\gamma).$$
Recall the Cartan involution $\oZ^\alpha\simeq\oZ^{\alpha,0}
\stackrel{\iota}{\longrightarrow}\oZ^{0,-\alpha}\simeq\oZ^\alpha$
of~\refss{version}. The following lemma is used in the next~\refs{invol}.

\lem{plusminus}
The following diagram commutes:
\eq{obvi}
\begin{CD}
(\BA^\beta\times\BA^\gamma)_{\on{disj}}\times_{\BA^\alpha}\oZ^\alpha @>\ff_{\beta,\gamma}>>
(\BA^\beta\times\BA^\gamma)_{\on{disj}}\times_{\BA^\beta\times\BA^\gamma}
(\oZ^\beta\times\oZ^\gamma)\\
@V\on{Id}\times\iota VV      @V\on{Id}\times\iota\times\iota VV\\
(\BA^\beta\times\BA^\gamma)_{\on{disj}}\times_{\BA^\alpha}\oZ^\alpha @>\ff_{\beta,\gamma}^+>>
(\BA^\beta\times\BA^\gamma)_{\on{disj}}\times_{\BA^\beta\times\BA^\gamma}
(\oZ^\beta\times\oZ^\gamma)
\end{CD}
\end{equation}
\elem

\prf Obvious. \epr

\sec{invol}{Cartan involution (Proof of \reft{main}(4))}

\ssec{involu}{Involution in coordinates}
Recall the modified coordinates $\by_{i,r}$ of~\refss{modi}, and
the Cartan involution $\iota:\ \oZ^\alpha\to\oZ^\alpha$ of~\refss{version}.

\prop{ioim}
The involution $\iota:\ \oZ^\alpha\to\oZ^\alpha$ in coordinates acts as follows:\\
$\iota:\ (w_{i,r},\by_{i,r})\mapsto(w_{i,r},\by_{i,r}^{-1})$
(equivalently, $(w_{i,r},y_{i,r})\mapsto(w_{i,r},y_{i,r}^{-1}
\prod_{j\ne i}Q_j^{-\langle\alpha_j,\check\alpha_i\rangle}(w_{i,r}))$).
\eprop

\prf
Recall that a $B$-structure on $\CF_G$ is encoded in
a collection
$\kappa_{\check\lambda}:\ \CL_{\check\lambda}\hookrightarrow\CV^{\check\lambda}_{\CF_G}$
of line subbundles satisfying the Pl\"ucker relations. Equivalently,
we can consider a collection $\kappa^*_{-w_0\check\lambda}:\
\CV^{\check\lambda}_{\CF_G}\twoheadrightarrow\ '\!\CL_{\check\lambda}$ of the quotient
line bundles satisfying the Pl\"ucker relations (we have
$'\!\CL_{\check\lambda}=\CL_{-w_0\check\lambda}^*$).
Similarly, a $B_-$-structure on $\CF_G$ is encoded in a collection of
line subbundles
$\kappa_{\check\lambda}^-:\ \CL_{\check\lambda}^-\hookrightarrow\CV^{\check\lambda}_{\CF_G}$
or equivalently, a collection of the quotient line bundles
$\kappa^{-*}_{-w_0\check\lambda}:\
\CV^{\check\lambda}_{\CF_G}\twoheadrightarrow\ '\!\CL_{\check\lambda}^-$.
Let $P_i$ (resp. $P_i^-$) be the $i$-type subminimal parabolic subgroup
containing $B$ (resp. $B_-$). Then a $B$-structure on $\CF_G$ induces a
$P_i$-structure on $\CF_G$ that gives rise to a 2-dimensional subbundle
$\CV_i\hookrightarrow\CV^{\check\omega_i}_{\CF_G}$ (associated to the 2-dimensional
subspace of invariants $V_{\check\omega_i}^{\on{Rad}P_i}\subset V_{\check\omega_i}$).
Similarly, a $B_-$-structure on $\CF_G$ induces a $P_i^-$-structure on $\CF_G$
that gives rise to a 2-dimensional quotient bundle
$\CV^{\check\omega_i}_{\CF_G}\twoheadrightarrow\ '\CV_i^-$.

We have the natural embedding $\CL_{\check\omega_i}\hookrightarrow\CV_i$ and
the natural projection $'\CV_i^-\twoheadrightarrow\ '\!\CL_{\check\omega_i}$.
We define the line bundle $\CM_i:=\CV_i/\CL_{\check\omega_i}$ so that we have
a short exact sequence $0\to\CL_{\check\omega_i}\to\CV_i\to\CM_i\to0$.
We define the line bundle $'\CM_i$ as the kernel of
$'\CV_i^-\twoheadrightarrow\ '\!\CL_{\check\omega_i}$ so that we have a short
exact sequence $0\to\ '\CM_i\to\ '\CV_i^-\to\ '\!\CL_{\check\omega_i}\to0$.
We also consider the composition
$\CL_{\check\omega_i}\hookrightarrow\CV^{\check\omega_i}_{\CF_G}\twoheadrightarrow\
'\CV_i^-$. We define $\CN_i$ as the cokernel of this composed map.
Note that generically over $\oZ^{\lambda,\mu}$ this composed map is an embedding
of the line bundle $\CL_{\check\omega_i}$, so that $\CN_i$ is a line bundle
as well, and we have a short exact sequence
$0\to\CL_{\check\omega_i}\to\ '\CV_i^-\to\CN_i\to0$. Given a general
$(\CF_G,\kappa_{\check\lambda},\kappa_{\check\lambda}^{-*})\in\oZ^{\lambda,\mu}$ such that
$\CN_i$ is a line bundle, we consider the following diagram:
\eq{comdi}
\begin{CD}
\CL_{\check\omega_i}@>>>\CV_i@>>>\CM_i\\
@|   @VVV  @V{\fQ}VV\\
\CL_{\check\omega_i}@>>> '\CV_i^-@>>>\CN_i
\end{CD}
\end{equation}
Here the rows are the above short exact sequences, the middle vertical map
is defined as the composition
$\CV_i\hookrightarrow\CV^{\check\omega_i}_{\CF_G}\twoheadrightarrow\ '\CV_i^-$,
and the right vertical map $\fQ$ is defined as follows.
Note that the trivialization of $\CF_G$ at $\infty\in C$ compatible with
the $B,B_-$-structures gives rise to
the trivializations of $\CL_{\check\omega_i},\CM_i,\CN_i$ at $\infty\in C$.
For degree reasons, $\CM_i$ is canonically isomorphic to
$\CO_C(\langle\lambda,-\check\omega+\check\alpha_i\rangle)$, and
$\CN_i$ is canonically isomorphic to
$\CO_C(\langle\alpha,\check\omega_i\rangle+
\langle\mu,-\check\omega_i+\check\alpha_i\rangle)$.
Finally $\fQ\in\Hom(\CM_i,\CN_i)=
\Gamma(C,\CO(\langle\alpha,2\check\omega_i-\check\alpha_i\rangle))$ is
defined as $\prod_{j\ne i}Q_j^{-\langle\alpha_j,\check\alpha_i\rangle}$.

\lem{210}
The diagram~\refe{comdi} commutes.
\elem
\prf
Straightforward.
\epr

Now given a general
$(\CF_G,\kappa_{\check\lambda},\kappa_{\check\lambda}^{-*})\in\oZ^{\lambda,\mu}\simeq
\oZ^\alpha$, the coordinates $w_{i,r}$ are nothing but the points
of $C$ where the line subbundles $\CL_{\check\omega_i}\hookrightarrow\ '\CV_i^-$
and $'\CM_i\hookrightarrow\ '\CV_i^-$ are not transversal. The trivialization
of $\CL_{\check\omega_i},\ '\CM_i$ at $\infty\in C$ gives rise to a canonical
trivialization of these line bundles restricted to $\BA^1=C\setminus\{\infty\}$.
Hence at a nontransversality point $w_{i,r}\in\BA^1$ we have two collinear
vectors in the fiber $'\CV_i^-|_{w_{i,r}}$, and the coordinate $y_{i,r}$ is
nothing but their ratio.

Since the Cartan involution $\oZ^\alpha\simeq\oZ^{\lambda,\mu}
\stackrel{\iota}{\longrightarrow}\oZ^{-\mu,-\lambda}\simeq\oZ^\alpha$ takes
$(\CF_G,\kappa_{\check\lambda},\kappa_{\check\lambda}^{-*})$ to
$(\CF_G,\kappa_{\check\lambda}^{-*},\kappa_{\check\lambda})$, and interchanges
the line bundles $\CL_{\check\omega_i},\CM_i$ with and without primes,
the proposition follows.
\epr

\sec{bound}{An equation of the boundary (Proof of \reft{main}(2))}

\ssec{y}{An equation in modified coordinates}
A regular function $F_\alpha$ on $Z^\alpha$ was constructed in~\cite[Section~4]{bf}
such that the divisor of $F_\alpha$ is the boundary divisor $\partial Z^\alpha$
(the multiplicities of various irreducible components of the boundary are 1
or $d$), see~\cite[Lemma~4.2]{bf}. Recall the modified coordinates $\by_{i,r}$
of~\refss{modi}.

\th{boun}
There is $c_\alpha\in\BC^*$ such that $c_\alpha F_\alpha=\prod_{i,r}\by_{i,r}^{d_i}=
\prod_{i,r}y_{i,r}^{d_i}\prod_{j\ne i}Q_j^{\alpha_j\cdot\alpha_i/2}(w_{i,r})$.
\eth

The rest of the section is devoted to the proof of the theorem.

\ssec{inv}{Invertible functions on zastava}
Let us denote the RHS of~\reft{boun} by $\fF_\alpha$. If we can prove that $Y_\alpha$ is a regular function on $Z^\alpha$ invertible on $\oZ^\alpha$ with
a correct order of vanishing at $\partial Z^\alpha$, then $\fF_\alpha/F_\alpha$
is a rational function on $Z^\alpha$ regular and nonvanishing at $\oZ^\alpha$ and
at the generic
points of the irreducible components of the divisor $\partial Z^\alpha$. Due to
normality of $Z^\alpha$~\cite[Corollary~2.10]{bf}, the ratio $\fF_\alpha/F_\alpha$
is a regular invertible function on $Z^\alpha$. Then according to the following
lemma, the ratio $\fF_\alpha/F_\alpha$ is a nonzero constant $c_\alpha$.

\lem{inve}
$\Gamma(Z^\alpha,\CO^*_{Z^\alpha})=\BC^*$.
\elem

\prf
Recall the factorization morphism $\pi_\alpha:\ Z^\alpha\to\BA^\alpha$.
Let $\Delta\subset\BA^\alpha$ be the diagonal divisor. For an off-diagonal
configuration $D\in\BA^\alpha$ the fiber $\pi_\alpha^{-1}(D)$ is isomorphic to
the $\langle\alpha,\check\rho\rangle$-dimensional affine space. Hence for
$f\in\Gamma(Z^\alpha,\CO^*_{Z^\alpha})$ the restriction of $f$ to any off-diagonal
fiber of $\pi_\alpha$ is constant. Hence $f=\bar{f}\circ\pi_\alpha$ for a certain
(invertible) function $\bar{f}$ on $\BA^\alpha$. Such $\bar{f}$ is necessarily
constant.
\epr

\ssec{a1}{Codimension one: $A_1$ and $A_1\times A_1$}
The order of vanishing of $\fF_\alpha$ at the generic points of
the irreducible components of $\partial Z^\alpha$ clearly coincides with that
of $F_\alpha$: see~\cite[Lemma~4.2]{bf}. We prove the regularity
of $\fF_\alpha$. Due to normality of $Z^\alpha$ it suffices to check
the regularity at the generic points of divisors $w_{i,r}=w_{j,s}$.
By the factorization property, it suffices to consider the case
$\alpha=\alpha_i+\alpha_j$. The case when $\alpha_i\cdot\alpha_j=0$ being
evident, we start with $i=j$. Then we can assume $\fg=\mathfrak{sl}_2$, so
that $Z^2_{\mathfrak{sl}_2}\simeq\BA^4=\{(Q_i=z^2+\sa_1z+\sa_2,R_i=\sfb_0z+\sfb_1)\}$.
We have $Q_i=(z-w_1)(z-w_2),\ R_i=(y_1(z-w_2)-y_2(z-w_1))/(w_1-w_2)$, so that
$y_1y_2$ is the resultant $R(Q_i,R_i)$: a regular function on
$Z^2_{\mathfrak{sl}_2}$, an equation of the boundary.

\ssec{a2}{Codimension one: $A_2$}
Next assume $i\ne j$, and $\alpha_i\cdot\alpha_j\ne0$, and $d_i=d_j$.
Then we can assume $\fg=\mathfrak{sl}_3$. Both fundamental representations
$V_{\check\omega_i},V_{\check\omega_j}$ of $\mathfrak{sl}_3$ are 3-dimensional.
The zastava space $Z^{\alpha_i+\alpha_j}_{\mathfrak{sl}_3}$
is formed by the polynomials with values in
$V_{\check\omega_i},V_{\check\omega_j}$ of the form $(z-w_i,y_i,u),(z-w_j,y_j,-u)$
such that $y_iy_j+(w_i-w_j)u=0$. We have
$\frac{y_i}{\sqrt{w_j-w_i}}\frac{y_j}{\sqrt{w_i-w_j}}=\sqrt{-1}u$: a regular
function on $Z^{\alpha_i+\alpha_j}_{\mathfrak{sl}_3}$, an equation of the boundary.

\ssec{b2}{Codimension one: $B_2$}
Next assume $i\ne j$, and $\alpha_i\cdot\alpha_j\ne0$, and $d_i=2,\ d_j=1$.
Then we can assume $\fg=\mathfrak{sp}_4$. The fundamental representation
$V_{\check\omega_i}$ (resp. $V_{\check\omega_j}$) is 4-dimensional
(resp. 5-dimensional). The Pl\"ucker coordinates for
$Z^{\alpha_i+\alpha_j}_{\mathfrak{sp}_4}$  are as follows:
$$\begin{array}{ccc}
c_{03}&\framebox{$b_{01}$}&z+A_1\\
\framebox{$b_{03}$}&c_{02}&\framebox{$z+A_2$}\\
c&\framebox{$b_{02}$}&b_{12}
\end{array}.$$
Here the boxed coordinates are the weight components of $V_{\check\omega_i}$, and the
remaining ones are the weight components of $V_{\check\omega_j}$. They are placed in the weight
lattice of $Sp(4)$. The origin of the weird notation is
in~\cite[Example~2.3.2]{fr}.
The Pl\"ucker equations are as follows.
First, we have the natural pairing $V_{\check\omega_j}\otimes V_{\check\omega_j}\to\BC$ coming from
$V_{\check\omega_j}\subset\Lambda^2V_{\check\omega_i}$, and $\Lambda^4V_{\check\omega_i}=\BC$. The $V_{\check\omega_j}$-valued polynomial must
be selforthogonal. The vanishing of the leading coefficient of the selfpairing
is $c=0$. The vanishing of the degree zero coefficient of the selfpairing is
$2c_{03}b_{12}-2c_{02}^2=0$. Second, we have the projection
$V_{\check\omega_i}\otimes V_{\check\omega_j}\to\Lambda^3V_{\check\omega_i}$ which must vanish on our polynomials. The vanishing
of the leading coefficient of the projection is
$c_{02}=b_{02},\ c_{03}=b_{03}$.
The vanishing of the degree zero coefficient of the projection is
$-c_{02}b_{02}+cA_2+b_{12}b_{03}=0,\ A_1b_{02}-A_2c_{02}-b_{01}b_{12}=0,\
b_{02}c_{03}-b_{01}c-b_{03}c_{02}=0,\ A_2c_{03}-A_1b_{03}+b_{01}c_{02}=0$.
Note that the
former quadratic equation is equivalent to the {\em first} quadratic Pl\"ucker
equation. All in all, we can take
$A_1,A_2,b_{12},b_{01},b_{02},b_{03}$ as independent coordinates, and we will
have three quadratic equations:
$b_{02}(A_1-A_2)=b_{01}b_{12},\ b_{03}(A_1-A_2)=b_{01}b_{02},\ b_{02}^2=b_{12}b_{03}$
(noncomplete intersection of three quadrics). To compare with the coordinates
of~\refss{igrek}: $w_i=-A_1,\ w_j=-A_2,\ y_i=b_{01},\ y_j=b_{12}$.

We have $(\frac{y_i}{\sqrt{w_j-w_i}})^2\frac{y_j}{w_i-w_j}=
-\frac{b_{01}^2b_{12}}{(w_i-w_j)^2}=-\frac{b_{01}b_{02}}{w_i-w_j}=-b_{03}$: a regular
function on $Z^{\alpha_i+\alpha_j}_{\mathfrak{sp}_4}$, an equation of the boundary.

\ssec{g2}{Codimension one: $G_2$}
Next assume $i\ne j$, and $\alpha_i\cdot\alpha_j\ne0$, and $d_i=3,\ d_j=1$.
Then $\fg$ is of type $G_2$. We have the regular functions $w_i,w_j,y_i,y_j$ on
$Z^{\alpha_i+\alpha_j}_\fg$. We have to show that
$(\frac{y_i}{\sqrt{w_j-w_i}})^3\frac{y_j}{\sqrt{(w_i-w_j)^3}}=
\sqrt{-1}\frac{y_i^3y_j}{(w_i-w_j)^3}$ is a regular function on
$Z^{\alpha_i+\alpha_j}_\fg$. According to the formulas of~\refss{modi}, the
Poisson bracket $\{y_i,y_j\}=-3\frac{y_iy_j}{w_i-w_j}$ is a regular function.
Furthermore, $\{y_i,\frac{y_iy_j}{w_i-w_j}\}=y_i\frac{\{y_i,y_j\}}{w_i-w_j}+
y_iy_j\{y_i,\frac{1}{w_i-w_j}\}=-3\frac{y_i^2y_j}{(w_i-w_j)^2}+
\frac{y_i^2y_j}{(w_i-w_j)^2}=-2\frac{y_i^2y_j}{(w_i-w_j)^2}$ is
a regular function. Finally,
$\{y_i,\frac{y_i^2y_j}{(w_i-w_j)^2}\}=-\frac{y_i^3y_j}{(w_i-w_j)^3}$ is a regular
function on $Z^{\alpha_i+\alpha_j}_\fg$.

\ssec{ani}{Invertibility}
The last thing to check is the invertibility of $\fF_\alpha$
on $\oZ^\alpha$. To this end recall the Cartan involution
$\iota:\ \oZ^\alpha\to\oZ^\alpha$ of~\refss{version} and note that
according to~\refp{ioim} we have
$\fF_\alpha\circ\iota=\fF_\alpha^{-1}$.

The theorem is proved. \qed

\rem{brave}
Here is an alternative way to prove the invertibility of
$\fF_\alpha$ on $\oZ^\alpha$: much shorter but less elementary.
According to~\cite[Proposition~4.4]{bf}, the weight of $F_\alpha$
with respect to the loop rotations in all the examples
of~\refss{a2},~\refss{b2},~\refss{g2} is equal to one.
If $\fF_\alpha$ were not regular on $Z^\alpha$, for certain $m>0$ the function
$(w_i-w_j)^m\fF_\alpha$ would be regular on $Z^\alpha$ and invertible on
$\oZ^\alpha$. Thus, the ratio $(w_i-w_j)^m\fF_\alpha/F_\alpha$
would be invertible on $Z^\alpha$
(since the numerator and denominator have the same order
of vanishing at the boundary $\partial Z^\alpha$) and hence constant
by~\refl{inve}. Since the weight of both $(w_i-w_j)$ and
$\fF_\alpha$ with respect to the loop rotations is also equal to one, we
conclude $m=0$: a contradiction with our assumption $m>0$. Hence we have
proved the regularity of $\fF_\alpha$ on $Z^\alpha$
and simultaneously the invertibility of $\fF_\alpha$
on $\oZ^\alpha$.
The general case is reduced to the above examples by factorization.
\erem

This completes the proof of~\reft{boun}. \qed

\sec{gaits}{An Ext calculation (Proof of \reft{main}(5))}

\ssec{pgl2}{$PGL_2$-bundles} A $PGL_2$-bundle with a flag on $C=\BP^1$ can be
viewed as a short exact sequence $0\to\CL\to\CV\to\CM\to0$ ($\CL$ and $\CM$
are the line bundles, and $\CV$ is a rank two vector bundle) modulo the
twistings by the line bundles. In particular, the line bundle
$\CM^{-1}\otimes\CL=\CHom(\CM,\CL)$ is well defined: this is nothing but the
induction of the Borel bundle to the Cartan bundle. We consider the moduli
stack $\CF_2$ of $PGL_2$-bundles with a flag on $C$ equipped with a
trivialization at $\infty\in C$ of the corresponding line bundle
$\CM^{-1}\otimes\CL$. The connected
components of $\CF_2$ are numbered by the integers $\deg\CM-\deg\CL$.
On a connected component $\CF_2^n$, we have a canonical isomorphism
$\CHom(\CM,\CL)=\CO_C(-n)$, and so
$\Ext^1(\CM,\CL)=H^1(C,\CO(-n))=H^0(C,\CO(n-2))^\vee$.
Thus we have a morphism $E:\ \CF_2^n\to H^0(C,\CO(n-2))^\vee$.

\ssec{mapext}{A map to Ext}
By Pl\"ucker, we may view a $B$-structure on $\CF_G$ as a collection of
line subbundles $\CL_{\check\lambda}\subset\CV^{\check\lambda}_{\CF_G}$ satisfying the
Pl\"ucker relations (here $\check\lambda$ runs through the cone of dominant
weights of $G$, and $\CV^{\check\lambda}_{\CF_G}$ is the vector bundle associated
to the irreducible $G$-module $V_{\check\lambda}$ with the highest weight
$\check\lambda$). For a $B$-structure coming from a point of $\oZ^{\lambda,\mu}$
we have $\deg\CL_{\check\lambda}=-\langle\lambda,\check\lambda\rangle$. The
trivialization at $\infty\in C$ extends to a canonical isomorphism
$\CL_{\check\lambda}=\CO_C(-\langle\lambda,\check\lambda\rangle)$.
Since the assignment $\check\lambda\mapsto\CL_{\check\lambda}$ is multiplicative
in $\check\lambda:\ \CL_{\check\mu+\check\nu}=\CL_{\check\mu}\otimes\CL_{\check\nu}$,
we can extend the notion of $\CL_{\check\lambda}$ for arbitrary weights
$\check\gamma\in X^*(T)$. We have a canonical isomorphism
$\CL_{\check\gamma}=\CO_C(-\langle\lambda,\check\gamma\rangle)$.

For $i\in I$ we have a morphism
$\frP_i:\ \oZ^{\lambda,\mu}\to\CF_2^{\langle\lambda,\check\alpha_i\rangle}$
defined as follows. Let $G\supset P_i\supset B$ be the $i$-type subminimal
parabolic
subgroup. We have the projection $P_i\to L_i$ to the corresponding Levi, also
we have the projection $L_i\to PGL_2$, and so the composed projection
$B\hookrightarrow P_i\twoheadrightarrow L_i\twoheadrightarrow PGL_2$.
Given a $B$-reduction of a $G$-bundle on $\CF_G$ trivialized at $\infty\in C$,
we consider the induced $PGL_2$-bundle. It comes
equipped with a flag and trivialization at $\infty\in C$. This is our desired
$\frP_i$ applied to $\CF_G$ with the $B$-structure (note that we have not used
the $B_-$-structure).

If we specialize to the case $(\lambda,\mu)=(\alpha,0)$, the assumption
$\mu=0$ guarantees that the $B_-$-structure further reduces to the
$U_-$-structure (where $U_-$ is the radical of $B_-$). Since any $U_-$-bundle
on $C=\BP^1$ is automatically trivial, the ambient $G$-bundle $\CF_G$ is
trivial too, and its trivialization at $\infty\in C$ extends to a canonical
trivialization over $C$. Thus we arrive at an identification
$\oZ^{\alpha,0}\simeq\oZ^\alpha$ with the usual zastava space, i.e. the moduli
space of degree $\alpha$ based maps from $(C,\infty)$ to $(\CB,B)$.
For $\phi\in\oZ^\alpha$ let us describe explicitly two particular
representatives (2-dimensional bundles with a flag) of $\frP_i(\phi)$.

We have a projection $p_i:\ \CB=G/B\to G/P_i=:\CB_i$. We define
$\CB^i:=\CB\times_{\CB_i}\CB$, and $p^i:\ \CB^i\to\CB$ (the first projection).
By construction, $p^i$ is a $\BP^1$-bundle over $\CB$ equipped with a
canonical (diagonal) section $\Delta^i:\ \CB\to\CB^i$.
We define $\CV'_i:=p^i_*\CO_{\CB^i}(\Delta^i)\supset\CL'_i:=p^i_*\CO_{\CB^i}=\CO_\CB$.
Thus we get a short exact sequence $0\to\CL'_i\to\CV'_i\to\CM'_i\to0$
trivialized at $B\in\CB$; here $\CM'_i=\CO_\CB(\check\alpha_i)$. Finally,
$\frP_i(\phi)=\{0\to\phi^*\CL'_i\to\phi^*\CV'_i\to\phi^*\CM'_i\to0\}$.

Alternatively, let $\CV_{\check\omega_i}$ be the trivial vector bundle over $\CB$
associated with the fundamental $G$-module $V_{\check\omega_i}$. It has a line
subbundle $\CL_i$: the fiber $\CL_i|_{B'}$ is the $B'$-highest line
$V_{\check\omega_i}^{\on{Rad}B'}$. If $P'_i$ is the $i$-type subminimal parabolic
containing $B'$, then the invariants $V_{\check\omega_i}^{\on{Rad}P'_i}$ are
2-dimensional (the highest and next highest lines), and as $B'$ varies in $\CB$,
we obtain a 2-dimensional subbundle $\CV_i\subset\CV_{\check\omega_i}$.
Thus we have a short exact sequence $0\to\CL_i\to\CV_i\to\CM_i\to0$ trivialized
at $B\in\CB$; here $\CL_i=\CO_\CB(-\check\omega_i)$, and
$\CM_i=\CO_\CB(-\check\omega_i+\check\alpha_i)$. Again, we have
$\frP_i(\phi)=\{0\to\phi^*\CL_i\to\phi^*\CV_i\to\phi^*\CM_i\to0\}$.

Finally, we define $\fE_i:\ \oZ^{\lambda,\mu}\to\Ext^1(\CO_C,\CL_{\check\alpha_i})=
\Ext^1(\CO_C,\CO_C(\langle\mu-\lambda,\check\alpha_i\rangle))$ as the
composition $\oZ^{\lambda,\mu}\simeq\oZ^{\lambda-\mu,0}\stackrel{\frP_i}{\longrightarrow}
\CF_2^{\langle\lambda-\mu,\check\alpha_i\rangle}\stackrel{E}{\longrightarrow}
\Ext^1(\CO_C,\CO_C(\langle\mu-\lambda,\check\alpha_i\rangle))$.

\ssec{memoir}{Recollections of~\cite{g}}
We recall some of the constructions of~\cite{g} in the particular case of
a curve of genus 0 (projective line $C$). In this case the canonical bundle
$\omega_C$ is isomorphic to $\CO_C(-2)$, and we choose a square root
$\omega_C^{1/2}\simeq\CO_C(-1)$.
Let $\Lambda=(\lambda_1,\ldots,\lambda_N)$ be an ordered collection of dominant
coweights. Let $\oA^\Lambda$ be the moduli space of ordered configurations of
distinct points $(z_1,\ldots,z_N\in\BA^1)$. Let
$\oA^{\alpha,\Lambda}\subset\BA^\alpha\times\oA^\Lambda$ be the open subspace formed
by the configurations of pairwise distinct points.
For $i\in I$ and $\ul{z}\in\oA^\Lambda$ we define a monic polynomial
$K_i(z):=\prod_{1\leq n\leq N}(z-z_n)^{\langle\lambda_n,\check\alpha_i\rangle}$.

Given a point $\unl{z}\in\oA^\Lambda$, we
consider a moduli stack $\fW_{\unl{z},\Lambda}$ classifying the following data:
(a) A $G$-bundle $\CF_G$ on $C$; (b) For each dominant weight $\check\lambda$
a nonzero map $\kappa^{\check\lambda}:\ \omega_C^{\langle\rho,\check\lambda\rangle}\to
\CV^{\check\lambda}_{\CF_G}$ having the poles of order exactly $\lambda_i$ at $z_i$,
and regular nonvanishing at $C\setminus\{\unl{z}\}$. Here
$\CV^{\check\lambda}_{\CF_G}$
is the vector bundle associated to $\CF_G$ and the irreducible $G$-module
$V^{\check\lambda}$, and $\omega_C^{\langle\rho,\check\lambda\rangle}$ stands for
$(\omega_C^{1/2})^{\otimes\langle2\rho,\check\lambda\rangle}$. The collection of maps
$\kappa^{\check\lambda}$ must satisfy the Pl\"ucker relations
(cf.~\cite[2.1,~2.6]{g}).

Alternatively, note that
$\kappa^{\check\lambda}|_{\omega_C^{\langle\rho,\check\lambda\rangle}(-\sum_{n=1}^N
\langle\lambda_n,\check\lambda\rangle\cdot z_n)}$ is a regular embedding,
and the image is a line subbundle
$\CL_{\check\lambda}\subset\CV^{\check\lambda}_{\CF_G}$.
So $\fW_{\unl{z},\Lambda}$ is the moduli stack of the collections of line
subbundles $\CL_{\check\lambda}\subset\CV^{\check\lambda}_{\CF_G}$ satisfying the
Pl\"ucker relations plus the identifications
$\CL_{\check\lambda}=\omega_C^{\langle\rho,\check\lambda\rangle}(-\sum_{n=1}^N
\langle\lambda_n,\check\lambda\rangle\cdot z_n)$.

Since the assignment $\check\lambda\mapsto\CL_{\check\lambda}$ is multiplicative
in $\check\lambda:\ \CL_{\check\mu+\check\nu}=\CL_{\check\mu}\otimes\CL_{\check\nu}$,
we can extend the notion of $\CL_{\check\lambda}$ for arbitrary weights
$\check\gamma\in X^*(T)$. The construction of~\refss{mapext} defines a morphism
$\fE_i:\ \fW_{\unl{z},\Lambda}\to\Ext^1(\CO_C,\CL_{\check\alpha_i})$.
The canonical embedding
$\CL_{\check\alpha_i}\hookrightarrow\omega_C^{\langle\rho,\check\alpha_i\rangle}=\omega_C$
gives rise to the projection
$\Ext^1(\CO_C,\CL_{\check\alpha_i})\to\Ext^1(\CO_C,\omega_C)=\BA^1$. Composing it
with $\fE_i$ we obtain a function $\chi_i:\ \fW_{\unl{z},\Lambda}\to\BA^1$.

Following~\cite[4.3]{g} we consider the moduli stack
$\oCZ^\alpha\to\fW_{\unl{z},\Lambda}$ classifying the same data as $\fW_{\unl{z},\Lambda}$
plus (a) a trivialization of $\CF_G$ at $\infty\in C$ such that the
$B$-structure (given by the collection $\{\kappa^{\check\lambda}\}$) at $\infty$
coincides with $B\subset G$; (b) an additional $B_-$-structure on $\CF_G$ of
degree $-2\rho-\alpha$ equal at $\infty\in C$ to $B_-\subset G$.
By an abuse of notation we preserve the
notation $\chi_i:\ \oCZ^\alpha\to\BA^1$ for the composition of
$\chi_i:\ \fW_{\unl{z},\Lambda}\to\BA^1$ and the projection
$\oCZ^\alpha\to\fW_{\unl{z},\Lambda}$. According to~\cite[4.5,~4.6]{g},
the stack $\oCZ^\alpha$ is actually a scheme; moreover, we have a canonical
isomorphism $\oCZ^\alpha=\oZ^{-2\rho,-2\rho-\alpha}=\oZ^\alpha$.
Thus we obtain a function
$\chi_i:\ \oZ^\alpha\to\BA^1$. If we allow $\unl{z}$ to vary in $\oA^\Lambda$,
we obtain the same named function $\chi_i$ on $\oZ^\alpha\times\oA^\Lambda$.

\th{gai}
The function $\chi_i$ on $\oZ^\alpha\times\oA^\Lambda$ in the coordinates
$(\unl{w},\unl{y},\unl{z})$ is given by
$$\chi_i(\unl{w},\unl{y},\unl{z})=\sum_{r=1}^{a_i}y_{i,r}^{-1}
\frac{\prod_{j\ne i}Q_j^{-\langle\alpha_j,\check\alpha_i\rangle}(w_{i,r})}
{Q'_i(w_{i,r})}K_i(w_{i,r})=\sum_{r=1}^{a_i}\by_{i,r}^{-1}
\frac{\prod_{j\ne i}Q_j^{-\langle\alpha_j,\check\alpha_i\rangle/2}(w_{i,r})}
{Q'_i(w_{i,r})}K_i(w_{i,r}).$$
\eth

\prf
Recall the involution $\iota:\ \oZ^\alpha\iso\oZ^\alpha$ of~\refss{version}.
We have to prove that $\chi_i\circ\iota(\unl{w},\unl{y},\unl{z})=
\sum_{r=1}^{a_i}y_{i,r}K_i(w_{i,r})/Q'_i(w_{i,r})$. Recall also the modified
coordinates $\fy_{i,r}:=\frac{y_{i,r}}{Q'_i(w_{i,r})}$ of~\refe{eta}.
Thus we have to prove
\eq{304}
\chi_i\circ\iota(\unl{w},\unl{\fy},\unl{z})=
\sum_{r=1}^{a_i}\fy_{i,r}K_i(w_{i,r})
\end{equation}
Recall the map
$\fE_i:\ \oZ^{\lambda,\mu}\to\Ext^1(\CO_C,\CL_{\check\alpha_i})$ of~\refss{mapext}.
We have to prove its factorization property, i.e. the commutativity of the
following diagram:
\eq{summa}
\begin{CD}
(\BA^\beta\times\BA^\gamma)_{\on{disj}}\times_{\BA^\alpha}\oZ^\alpha @>\ff_{\beta,\gamma}^+>>
(\BA^\beta\times\BA^\gamma)_{\on{disj}}\times_{\BA^\beta\times\BA^\gamma}
(\oZ^\beta\times\oZ^\gamma)\\
@V\fE_i VV      @V\fE_i\times\fE_i VV\\
\Ext^1(\CO_C,\CL_{\check\alpha_i})
@<+<<
\Ext^1(\CO_C,\CL_{\check\alpha_i})\times\Ext^1(\CO_C,\CL_{\check\alpha_i})
\end{CD}
\end{equation}

Unraveling the definition of $\fE_i$ and using the compatibility of
factorizations~\refe{cD} we reduce the problem to $G=SL_2^i$. This problem
is formulated as follows. For a point $x\in C$ and a line bundle
$\CE$ on $C$ we have a natural map
$\varphi_x:\ (\CE\otimes_{\CO_C}\CK_x)/(\CE\otimes_{\CO_C}\CO_x)\to H^1(C,\CE)$
arising from
the identification $(\CE\otimes_{\CO_C}\CK_x)/(\CE\otimes_{\CO_C}\CO_x)=
j_*\CE/\CE$ and the boundary map in the long exact cohomology sequence
coming from $0\to\CE\to j_*\CE\to j_*\CE/\CE\to0$ (here $j$ is the open
embedding $C\setminus\{x\}\hookrightarrow C$). Given a short exact sequence
$0\to\CL\to\CV\to\CM\to0$ as in~\refss{pgl2}, its local splittings form
a torsor $\CT$ over $\CE=\CHom(\CM,\CL)$ with the class
$e(\CT)\in H^1(C,\CE)$. Given a generic splitting of this
exact sequence we obtain a generic section $s$ of $\CT$. Choosing a
local splitting around $x$ we obtain $s_x\in\CE\otimes_{\CO_C}\CK_x$ whose
principal part in $(\CE\otimes_{\CO_C}\CK_x)/(\CE\otimes_{\CO_C}\CO_x)$ is well
defined, i.e. independent of the choice of a local splitting.
Then clearly $e(\CT)=\sum_{x\in C}\varphi_x(s_x)$.

This completes the proof of~\reft{gai}.
\epr

\ssec{Dmod}{More recollections (Proof of~\reft{gaits-to-witten})}
Recall the sequence of morphisms
$\oZ^\alpha\simeq\oZ^{-2\rho,-2\rho-\alpha}=\oCZ^\alpha\to\fW_{\unl{z},\Lambda}\to
\on{Bun}_G(C)$ of~\refss{memoir}. The line bundle $\CP$ on $\oZ^\alpha$ is
defined as the inverse image of the determinant line bundle on $\on{Bun}_G(C)$,
cf.~\cite[2.2]{g}. In fact, $\CP$ is the restriction of the same named
line bundle on $Z^\alpha$ with the canonical section $F_\alpha$,
see~\cite[4.9]{bf}. Hence $F_\alpha$ gives rise to a canonical trivialization
of $\CP$ on $\oZ^\alpha$. Given a level $\varkappa\in\BC$, Gaitsgory constructs
a certain $\CP^\varkappa$-twisted $D$-module $\CF_{\unl{z},\Lambda}^\varkappa$ on
$\oZ^\alpha\times\oA^\Lambda$ (as a lift from
$\fW_{\unl{z},\Lambda}\times\oA^\Lambda$)~\cite[2.7]{g}. It is smooth
of rank 1 on $\oZ^\alpha\times\oA^\Lambda$ but has irregular singularities at
$\partial Z^\alpha\times\oA^\Lambda$.
In case $\varkappa$ is irrational, $\CF_{\unl{z},\Lambda}^\varkappa$ is clean.
The trivialization of $\CP$ on $\oZ^\alpha$ gives rise to the identification
of $\CP^\varkappa$-twisted $D$-modules with the usual $D$-modules, and then
the corresponding $D$-module $\CF_{\unl{z},\Lambda}^{\varkappa,\on{triv}}$ on
$\oZ^\alpha\times\oA^\Lambda$ is generated by the function
$F_\alpha^\varkappa\cdot\exp(\chi_i)\cdot
\prod_{1\leq m<n\leq N}(z_m-z_n)^{\varkappa\lambda_m\cdot\lambda_n}$.

According to~\cite[Lemma~4.9]{g} there is a canonical isomorphism
$\CP=\pi_\alpha^*\CP_{\BA^\alpha}$ for a certain line bundle $\CP_{\BA^\alpha}$ on
$\BA^\alpha$~\cite[3.2]{g}. For irrational $\varkappa$~\cite[Theorem~6.2]{g}
identifies
$\pi_{\alpha!}\CF_{\unl{z},\Lambda}^\varkappa\iso\pi_{\alpha*}\CF_{\unl{z},\Lambda}^\varkappa$
as the minimal extension $\CL_{\unl{z},\Lambda}^\varkappa$ of a smooth rank
1 $\CP_{\BA^\alpha}^\varkappa$-twisted $D$-module from the open diagonal stratum of
$\BA^\alpha\times\oA^\Lambda$. The trivialization~\cite[3.12]{g} of
$\CP_{\BA^\alpha}$ on $\BA^\alpha\times\oA^\Lambda$ gives rise to the identification
of $\CP_{\BA^\alpha}^\varkappa$-twisted $D$-modules with the usual $D$-modules,
and then the corresponding $D$-module $\CL_{\unl{z},\Lambda}^{\varkappa,\on{triv}}$
is the minimal extension of the $D$-module on the open stratum generated by
the function
$$\prod_{(i,r)\ne(j,s)}(w_{i,r}-w_{j,s})^{\varkappa\alpha_i\cdot\alpha_j/2}\times
\prod_{(i,r),1\leq n\leq N}(z_n-w_{i,r})^{-\varkappa\alpha_i\cdot\lambda_n}\times
\prod_{1\leq m<n\leq N}(z_m-z_n)^{\varkappa\lambda_m\cdot\lambda_n}.$$

More generally, we consider the $D$-module $M^{\varkappa,\alpha,\Lambda}$ on
$\fh^\vee\times\oZ^\alpha\times\oA^\Lambda$ generated by the function
$\prod_{1\leq n\leq N}\exp(\langle\lambda_n,\varkappa h^*\rangle z_n)\cdot
\prod_{(i,r)}\exp(-\langle\alpha_i,\varkappa h^*\rangle w_{i,r})\cdot
F_\alpha^\varkappa\cdot\exp(\chi_i)\cdot
\prod_{1\leq m<n\leq N}(z_m-z_n)^{\varkappa\lambda_m\cdot\lambda_n}$.
Then for irrational $\varkappa$ the direct image
$\pi_{\alpha!}M^{\varkappa,\alpha,\Lambda}\iso\pi_{\alpha*}M^{\varkappa,\alpha,\Lambda}$
is isomorphic to the minimal extension
$\CM^{\varkappa,\on{triv}}_{\fh^\vee,\unl{z},\Lambda}$
of the $D$-module on the open diagonal
stratum of $\fh^\vee\times\BA^\alpha\times\oA^\Lambda$ generated by the function
$$\prod_{1\leq n\leq N}\exp(\langle\lambda_n,\varkappa h^*\rangle z_n)\times
\prod_{(i,r)}\exp(-\langle\alpha_i,\varkappa h^*\rangle w_{i,r})\times$$
$$\times\prod_{(i,r)\ne(j,s)}(w_{i,r}-w_{j,s})^{\varkappa\alpha_i\cdot\alpha_j/2}\times
\prod_{(i,r),1\leq n\leq N}(z_n-w_{i,r})^{-\varkappa\alpha_i\cdot\lambda_n}\times
\prod_{1\leq m<n\leq N}(z_m-z_n)^{\varkappa\lambda_m\cdot\lambda_n}.
$$
In effect, the isomorphism
$\pi_{\alpha!}M^{\varkappa,\alpha,\Lambda}\iso\pi_{\alpha*}M^{\varkappa,\alpha,\Lambda}\iso
\CM^{\varkappa,\on{triv}}_{\fh^\vee,\unl{z},\Lambda}$ follows from the isomorphism
$\pi_{\alpha!}\CF_{\unl{z},\Lambda}^{\varkappa,\on{triv}}\iso
\pi_{\alpha*}\CF_{\unl{z},\Lambda}^{\varkappa,\on{triv}}\iso
\CL_{\unl{z},\Lambda}^{\varkappa,\on{triv}}$ and the projection formula.
The latter isomorphism is proved in~\cite{g} for any {\em fixed} value of
$\unl{z}$. To prove it for variable $\unl{z}$ it remains to identify the
monodromy of the one-dimensional local system
$\pi_{\alpha!}M^{\varkappa,\alpha,\Lambda}$ on the open diagonal stratum.
This follows from the computation of the~\refp{ide} below.

\ssec{mas}{The master function and the Gaiotto-Witten superpotential}
The (multivalued) Master function~\cite[Section~3]{fmtv} on
$\fh^\vee\times\oA^{\alpha,\Lambda}$ is defined as follows:
\begin{multline}
\label{mast}
\Phi(h^*,\ul{w},\ul{z}):=\sum_{1\leq n\leq N}\langle\lambda_n,h^*\rangle z_n-
\sum_{(i,r)}\langle\alpha_i,h^*\rangle w_{i,r}+
\sum_{(i,r)\ne(j,s)}\frac{\alpha_i\cdot\alpha_j}{2}\log(w_{i,r}-w_{j,s})-\\
-\sum_{(i,r),1\leq n\leq N}\alpha_i\cdot\lambda_n\log(z_n-w_{i,r})+
\sum_{1\leq m<n\leq N}\lambda_m\cdot\lambda_n\log(z_m-z_n)
\end{multline}

We define an open subvariety $\oZ^{\alpha,\Lambda}\subset\oZ^\alpha\times\oA^\Lambda$
as the preimage of $\oA^{\alpha,\Lambda}$ under the factorization morphism
$\pi_\alpha:\ \oZ^\alpha\times\oA^\Lambda\to\BA^\alpha\times\oA^\Lambda$.
Recall the logarithmic coordinates $\sy_{i,r}$ of~\refd{loga}.
The multivalued superpotential $\CW^{\Lambda,\alpha}$ on
$\fh^\vee\times\oZ^{\alpha,\Lambda}$
is defined as follows (cf.~\cite{gw}):
\begin{multline}
\label{supe}
\CW^{\Lambda,\alpha}(h^*,\ul{w},\ul\sy,\ul{z}):=
\sum_{1\leq n\leq N}\langle\lambda_n,h^*\rangle z_n-
\sum_{(i,r)}\langle\alpha_i,h^*\rangle w_{i,r}+\sum_{(i,r)}d_i\sy_{i,r}+\\
+\sum_{(i,r)}\exp(-\sy_{i,r})
\frac{\prod_{j\ne i}Q_j^{-\langle\alpha_j,\check\alpha_i\rangle/2}(w_{i,r})}
{Q'_i(w_{i,r})}K_i(w_{i,r})+
\sum_{1\leq m<n\leq N}\lambda_m\cdot\lambda_n\log(z_m-z_n)
\end{multline}

\prop{ide}
a) The restriction of $\CW^{\Lambda,\alpha}$
to each fiber of the factorization projection
$\fh^\vee\times\oZ^{\alpha,\Lambda}\to\fh^\vee\times\oA^{\alpha,\Lambda}$ has a
unique singular point (with all the derivatives vanishing).

b) For any $h^*\in\fh^\vee,\ \ul{z}\in\oA^\Lambda$, the resulting section
$s_{h^*,\ul{z}}:\ \oA^\alpha\hookrightarrow\oZ^\alpha$ is Lagrangian.

c) The restriction of $\CW^{\Lambda,\alpha}$ to the section in a) equals the
Master function $\Phi$.
\eprop

\prf
Straightforward.
\epr

\sec{obso}{Appendix}
In this appendix we give another (elementary) derivation of a particular case
of~\reft{gai} for $G=SL(2)$.

\ssec{another}{Zastava for $SL(2)$}
For $G=SL(2)$ the coroot lattice is just $\BZ$, and given $a\in\BN$
the moduli space of based maps $\oZ^a$ is identified with the moduli
space of extensions $0\to\CO_C(-a)\to\CO_C\oplus\CO_C\to\CO_C(a)\to0$
trivialized at $\infty\in C$. So we have a map $\fE:\ \oZ^a\to
\Ext^1(\CO_C(a),\CO_C(-a))=\Gamma(C,\CO_C(2a-2))^\vee$.

\prop{formula'}
For a polynomial $K\in
H^0(C,\CO_C(2a-2))$ we have
$\displaystyle{\langle K,\fE\rangle=\sum_{r=1}^ay_r^{-1}
\frac{K(w_r)}{Q'(w_r)}.}$
\eprop

\prf
We denote $H^0(C,\CO_C(1))^\vee$ by $V$
(a 2-dimensional vector space with a base formed by the highest vector $x$ and
the lowest vector $t$). We have
$\Ext^1(\CO(a),\CO(-a))=\Gamma(C,\CO(2a-2))^\vee=\Sym^{2a-2}V$.
We will write down an element of $\Sym^{2a-2}V$ in the basis of products of
divided powers of
$x,t:\ c_0x^{(2n-2)}+\ldots+c_kx^{(2a-2-k)}t^{(k)}+\ldots+c_{2n-2}t^{(2a-2)}$.

For a point $\phi\in\oZ^a$, the first map in the corresponding exact sequence
$0\to\CO_C(-a)\to\CO_Cx\oplus\CO_Ct\to\CO_C(a)\to0$
is given by a pair of polynomials $(Q,R)$, and the second one is given by
$(-R,Q)$. In the corresponding long exact sequence
$0=H^0(\CO_C(-a))\to H^0(\CO_C\oplus\CO_C)\to H^0(\CO_C(a))\to
H^1(\CO_C(-a))\to\ldots$
the boundary map is given by the cup product with our desired $\on{Ext}^1$-class
$\fE(\phi)$ in $\on{Sym}^{2a-2}V$. Note that $H^0(\CO_C(a))=\on{Sym}^aV^\vee$, and
$H^1(\CO_C(-a))=\on{Sym}^{a-2}V$. So the boundary map is the contraction
$\on{Sym}^aV^\vee\to\on{Sym}^{a-2}V$ with the desired class
$\fE(\phi)\in\on{Sym}^{2a-2}V$.
Since the composition $H^0(\CO_C\oplus\CO_C)\to H^0(\CO_C(a))\to H^1(\CO(-a))$
is 0, and the first map is given by $(-R,Q)$, we conclude that the
contraction of $\fE(\phi)$ and $Q$ equals 0,
as well as the contraction of $\fE(\phi)$ and $R$.
This is a system of linear equations on $\fE(\phi)$ which defines it up to
proportionality. To write down the formula for contraction, we think of
$Q,R$ as of differential operators
$Q=\partial_x^a+\ldots+\sa_k\partial_x^{a-k}\partial_t^k+\ldots+\sa_a
\partial_t^a,\
R=\sfb_0\partial_x^{a-1}+\ldots+\sfb_k\partial_x^{a-1-k}\partial_t^k+\ldots+
\sfb_a\partial_t^{a-1}$, and then the contraction is nothing but the application
of differential operators $Q,R$ to the polynomial
$\fE:=c_0x^{(2a-2)}+\ldots+c_kx^{(2a-2-k)}t^{(k)}+\ldots+c_{2a-2}t^{(2a-2)}$.

Note that the matrix of this system of equations is (up to proportionality)
exactly the Sylvester matrix
$S=\left(
\begin{array}{ccccccc}
1&\sa_1&\ldots\sa_a&0&\ldots&0\\
\vdots&\ddots&\ddots&\ddots&\ddots&\ddots&\vdots\\
0&\ldots&0&1&\sa_1&\ldots&\sa_a\\
\sfb_0&\sfb_1&\ldots&\sfb_{a-1}&0&\ldots&0\\
\vdots&\ddots&\ddots&\ddots&\ddots&\ddots&\vdots\\
0&\ldots&0&\sfb_0&b_1&\ldots&\sfb_{a-1}
\end{array}
\right)$
with the middle row (the first one with $\sfb$'s) removed.
Solving it via the Cramer rule we obtain $c_k=(-1)^k\det S^{-1}$ times
the $(2a-2)\times(2a-2)$-minor of the Sylvester matrix obtained by removing
the middle row and the $k$-th column. Note also that the resultant $\sR(Q,R)$
is nothing but $\det S$, and $\sR(Q,R)\ne0$ under our assumptions:
$(\CO_C\oplus\CO_C)/\CO_C(-a)$ torsionless.

Equivalently, if we think of $Q,R$ as of two (relatively prime)
polynomials in $z=\partial_x/\partial_t$ (as in~\refss{igrek}), then the
equation $RD-QF=1$ has a unique solution such that $D$ is a polynomial in
$z$ of degree $a-1$, and $F$ is a polynomial in $z$ of degree $a-2$.
The principal part at $\infty\in C$ of the ratio $\frac{D(z)}{Q(z)}$
is nothing but
$\frac{c_0}{z}+\frac{c_1}{z^2}+\ldots+\frac{c_{2a-2}}{z^{2a-1}}+\ldots$ ($c_k$ from
the previous paragraph).
By the Lagrange interpolation we find
$c_k=\sum_{r=1}^a\frac{w_r^kD(w_r)}{Q'(w_r)}=
\sum_{r=1}^a\frac{w_r^kR(w_r)^{-1}}{Q'(w_r)}=
\sum_{r=1}^a\frac{w_r^ky_r^{-1}}{Q'(w_r)}$.
The desired formula for $\langle K,\fE\rangle$ follows.
\epr

\rem{kro}
We keep the notations introduced in the proof of~\refp{formula'}.
Let us define $\tilde{c}_0,\ldots,\tilde{c}_{2a-2}$ by
$\frac{R(z)}{Q(z)}=\frac{\tilde{c}_0}{z}+\frac{\tilde{c}_1}{z^2}+\ldots+
\frac{\tilde{c}_{2a-2}}{z^{2a-1}}+\ldots$ Then by the Lagrange interpolation
$\tilde{c}_k=\sum_{r=1}^a\frac{w_r^ky_r}{Q'(w_r)}$.
According to L.~Kronecker~\cite{k}, the resultant $\sR(Q,R)=\det\tilde{L}$
where $\tilde{L}$ is a Hankel matrix $\tilde{L}:=\left(
\begin{array}{ccccc}
\tilde{c}_0&\tilde{c}_1&\tilde{c}_2&\ldots&\tilde{c}_{a-1}\\
\tilde{c}_1&\tilde{c}_2&\tilde{c}_3&\ldots&\tilde{c}_a\\
\tilde{c}_2&\tilde{c}_3&\tilde{c}_4&\ldots&\tilde{c}_{a+1}\\
\vdots&\vdots&\vdots&\ddots&\vdots\\
\tilde{c}_{a-1}&\tilde{c}_a&\tilde{c}_{a+1}&\ldots&\tilde{c}_{2a-2}
\end{array}
\right)$. We obtain $\sR(Q,R)=\sR(D,R)^{-1}=\det L^{-1}$ where
$L$ is a Hankel matrix $L:=\left(
\begin{array}{ccccc}
c_0&c_1&c_2&\ldots&c_{a-1}\\
c_1&c_2&c_3&\ldots&c_a\\
c_2&c_3&c_4&\ldots&c_{a+1}\\
\vdots&\vdots&\vdots&\ddots&\vdots\\
c_{a-1}&c_a&c_{a+1}&\ldots&c_{2a-2}
\end{array}
\right)$. This identity $\sR(Q,R)=\det L^{-1}$ was independently obtained by
A.~Uteshev (private
communication). Note that the equation $\det L=0$ is the equation of the
locus in $\Ext^1(\CO_C(a),\CO_C(-a))$ formed by the extensions with the middle
term a {\em nontrivial} 2-dimensional vector bundle on $C$~\cite{cs}.
\erem

\bigskip

\footnotesize{ {\bf A.B.}: Department of Mathematics, Brown University,
151 Thayer st, Providence RI 02912, USA;\\ 
{\tt braval@math.brown.edu}}

\footnotesize{ {\bf G.D.}: Department of Mathematics, Columbia University,
New York, NY 10027, USA;\\
{\tt galdobr@gmail.com}}

\footnotesize{ {\bf M.F.}: National Research University Higher
School of Economics, Math. Dept.,\\ 20 Myasnitskaya st, Moscow 101000
Russia; IITP\\ 
{\tt fnklberg@gmail.com}}

\end{document}